\documentclass{article}

\usepackage{latexsym,amssymb,amscd,epic,eepic,epsfig,graphics,pstricks,verbatim
,fancyhdr}
\usepackage{amsbsy,amsmath,amsfonts}
\usepackage{amsthm}
\usepackage[francais,english]{babel}
\usepackage{graphicx}
\usepackage[all]{xy}
\usepackage{array}
\usepackage{enumerate}
\usepackage{amscd}

\newtheorem{theoremIntro}{Theorem}
\newtheorem{corollaryIntro}[theoremIntro]{Corollary}

\newtheorem{theorem}{Theorem}[section]
\newtheorem{lemma}[theorem]{Lemma}
\newtheorem{corollary}[theorem]{Corollary}
\newtheorem{proposition}[theorem]{Proposition}

\newtheorem{definition}[theorem]{Definition}

\newtheorem{remark}[theorem]{Remark}

\newtheorem*{NB*}{Nota Bene}

\newenvironment{sproof}[1]
{\begin{proof}[#1]} {\end{proof}}

\newcommand{\ov}[1]{\overline{#1}}
\newcommand{\N}{\mathbb N}
\newcommand{\Nast}{\mathbb{N}^{\ast}}
\newcommand{\Zast}{\mathbb{Z}^{\ast}}

\newcommand{\Za}{{\mathbb Z}^{\ast}}

\newcommand{\Z}{\mathbb Z}
\newcommand{\R}{\mathbb R}
\newcommand{\Q}{\mathbb Q}

\newcommand{\Free}{\mathbb{F}}
\newcommand{\F}{{\mathbb F}_2}
\newcommand{\Fm}{{\mathbb F}_m}

\newcommand{\G}{{\mathcal G}_2}
\newcommand{\Gm}{{\mathcal G}_m}

\newcommand{\BS}{\overline{BS}(m,\xi)}

\newcommand{\bBS}[2]{\overline{BS}(#1,#2)}

\newcommand{\BBS}{\overline{BS}}
\newcommand{\m}{m}
\newcommand{\n}{n}

\newcommand{\Pres}[2]{\left\langle{#1}\ \big\vert\ {#2}\right\rangle}
\newcommand{\eg}[1]{\underset{#1}{=}}

\title{Limits of Baumslag-Solitar groups }\label{ChAnnexeBS}
\author{Luc Guyot\footnote{Supported by the Swiss National Science Foundation,
    No.~PP002-68627.} \ and Yves Stalder\footnote{Supported by the Swiss National Science Foundation,
grant  20-109130, and the Universit\'e de Neuch\^atel, where he was employed during the research presented
in this article.}}
\date{March 21, 2008}

\begin{document}

\maketitle



\begin{abstract}
We give a parametrization by $m$-adic integers of the limits of
Baumslag-Solitar groups (marked with a canonical set of generators).
It is shown to be continuous and injective on the invertible
$m$-adic integers. We show that all such limits are extensions of a
free group by a lamplighter group and all but possibly one are not
finitely presented. Finally, we give presentations related to
natural actions on trees.
\end{abstract}


\section*{Introduction} \label{SecIntroAnnexeA}

The set $\mathcal G_k$ of marked groups on $k$ generators (see
{Section \ref{DefPrelim}} for definitions) has a natural topology,
which turns it into a metrizable, compact and totally disconnected
space. This topology, introduced explicitly by Grigorchuk
\cite{Gri84}\footnote{Grigorchuk and Zuk also defined spaces of marked graphs \cite{GZ99}.}, corresponds to an earlier and more general
construction by Chabauty \cite{Chab50}.

Grigorchuk topology was introduced in the context of growth of
finitely generated groups\footnote{Gromov used a very similar topology while discussing further its celebrated
polynomial growth theorem \cite[Final Remarks]{Grom81}.}. Grigorchuk
constructed a set of $3$-generated groups in $\mathcal{G}_4$ that is
homeomorphic to the Cantor set and he fruitfully studied the
neighborhoods of the \it first (intermediate growth) Grigorchuk group
\rm inside this set \cite{Gri84}. The Grigorchuk topology has
various connections with problems in group theory:
\begin{itemize}
    \item Using the Cantor of Grigorchuk's groups, Stepin proved the existence of uncountably many amenable but
    non-elementary amenable groups without growth estimates \cite{Ste84}.
    \item Shalom proved that every finitely generated Kazhdan
    group is a quotient of some finitely \textit{presented} Kazhdan group by
    showing that Kazhdan's property (T) defines an open subset of the
    space of marked groups \cite{Sha00}.
    \item Champetier showed that the quotient of the space
    of marked groups on $k$ generators by the group isomorphism relation is not
    a standard Borel space. He also proved that the closure of
    non-elementary hyperbolic groups contains a $G_\delta$ dense
    subset consisting of infinite torsion groups \cite{Cham00}.
    \item Champetier and Guirardel characterized \emph{limit
    groups} of Sela (this class coincides with fully residually free groups and
    with groups having the same universal theory as a free group) as limits of
    free groups. They also related Grigorchuk topology to universal theory of
    groups and ultraproducts \cite{ChGu05}.
    \item Nekrashevych constructed a minimal Cantor set
    in $\mathcal{G}_3$ \cite{Nek07} and proved that it contains a group with
    non-uniform exponential growth locally isomorphic to the
    iterated monodromy group of $z^2+i$ \cite{Nek10}.
\end{itemize}

We are interested in the closure of the set of Baumslag-Solitar groups
equipped with their standard marking and we study individual elements of this closure for their own
right. Let us recall that Baumslag-Solitar groups are defined by
\[
BS(\m,\n) = \Pres{a,b}{ab^\m a^{-1} = b^\n} \ \text{ for } \m,\n \in \Z\setminus\{0\}  \ .
\]
These groups have been introduced in \cite{BS62}, in order to give
the first examples of non-Hopfian one-relator groups. Since then,
they received much attention and served as test cases for several
problems in combinatorial and geometric group theory. They bear
further pathologies in various fields such as topology, geometry and
algebraic geometry over groups. Like $BS(2,3)$, many of them:
\begin{itemize}
\item do not embed into the fundamental group of any sufficiently large, irreducible, compact and connected $3$-manifold
\cite{JS79};
\item have no proper action on a CAT(0) cube complex
\cite{Hag07};
\item are not equationally Noetherian \cite{BMR99}.
\end{itemize}
Baumslag-Solitar groups were classified up to group isomorphism by
Moldavanski\u\i \ \cite{Mol91} and up to quasi-isometry by Farb and
Mosher \cite{FM98} and Whyte \cite{Why01}.

This paper is a continuation of \cite{Sta06a}, where the
second-named author characterized convergent sequences among
Baumslag-Solitar groups marked by generators $a$ and $b$, provided
the two parameters $\m$ and $\n$ are coprime. Given any non-zero
integer $m$, let us recall that the topological ring $\Z_\m$ of
$m$-adic integers is the projective limit of $\Z/\m^h\Z$ for
$h\in\N^*$ (see Section \ref{mAdicSct} or \cite[Ch. II.\S.10]{HR63}
for more details). Given any $\xi$ in $\Z_\m$, using results of
\cite{Sta06a}, we define:
\[
\bBS{\m}{\xi} = \lim_{n \to \infty} BS(\m,\xi_n)
\]
where $(\xi_n)_n$ is a sequence of rational integers such that
$\xi_n$ tends to $\xi$ in $\Z_\m$ and $|\xi_n|$ tends to infinity as
$n$ goes to infinity. This yields a parametrization
$\overline{BS}_\m : \Z_\m \rightarrow \mathcal G_2$ defined by $\xi
\mapsto \bBS{\m}{\xi}$. In Section 2, we prove the following:
\begin{theoremIntro}[Theorem \ref{ThmBSm} and Corollary \ref{CorBSm}]
Let $\m \in \Z \setminus \{0\}$ and  let $\xi,\eta \in \Z_{\m}$. The
equality of marked groups
$\overline{BS}(\m,\xi)=\overline{BS}(\m,\eta)$ holds if and only if
there is some $d\in\Z \setminus \{0\}$ such that
$gcd(\xi,\m)=gcd(\eta,\m)=d$ and the images of $\xi/d$ and $\eta/d$
in $\Z_{\m/d}$ are equal. In particular, the map $\overline{BS}_\m$
is injective on the set of invertible $\m$-adic integers.
\end{theoremIntro}
\begin{theoremIntro}[Corollary \ref{Embedding}]
The map $\overline{BS}_\m $ is continuous.
\end{theoremIntro}
These theorems enable us to describe (in Corollary \ref{Bdy}) the set
of groups $\bBS{\m}{\xi}$ for which $\xi$ is invertible as the boundary of a
particular set of Baumslag-Solitar groups in $\G$.

It is well-known that a Baumslag-Solitar group acts on its
Bass-Serre tree by automorphisms and on $\Q$ by affine
transformations. The first action is faithful whereas the second is
not in general. Section \ref{Actions} is devoted to actions
and structure of the groups $\bBS{\m}{\xi}$.
The groups $\bBS{\m}{\xi}$ mainly act by automorphisms on a tree which
is in some sense a ``limit'' of Bass-Serre trees (Theorem
\ref{TreeAction}). These considerations imply:
\begin{theoremIntro}[Theorem \ref{Extension}]\label{StructIntro} For any
$m\in\Z\setminus \{0\}$ and $\xi\in\Z_m$,
there exists an exact sequence $1 \rightarrow \Free \rightarrow
\bBS{\m}{\xi} \rightarrow \Z\wr\Z \rightarrow 1$, where $\Free$ is a
free group.
\end{theoremIntro}
Note that some properties obviously hold for the groups
$\bBS{\m}{\xi}$, since these properties define closed sets in the
space of marked groups: the groups $\bBS{\m}{\xi}$ are torsion free
and centerless, their subgroup generated by $a$ and $bab^{-1}$ is a
non-abelian free group for $\vert m \vert \ge 2$, and they are non-
Kazhdan. We now present some consequences of our Theorem
\ref{StructIntro}, which are not obvious in the former sense.
\begin{corollaryIntro}[Corollary \ref{Haagerup}]\label{HaIntro}
The limits $\bBS{\m}{\xi}$:
\begin{enumerate}
 \item have the Haagerup property;
 \item are residually solvable.
\end{enumerate}
\end{corollaryIntro}

In the last section of the paper, we discuss presentations of
the limits $\bBS{\m}{\xi}$. We prove in particular:
\begin{theoremIntro}[Theorem \ref{infPres}]
 For any $\m\in\Z \setminus \{0\}$ and $\xi\in\Z_\m \setminus \m\Z_\m$, the group $\ov{BS}(\m,\xi)$ is not finitely presented.
\end{theoremIntro}
We then exhibit with Theorem \ref{ThmPres1} presentations of the
limits, using again their actions by tree automorphisms.

\paragraph{Acknowledgements.} We would like to thank Thierry Coulbois for
having pointed out Remark \ref{ultra} to us. We also thank our
advisors, Goulnara Arzhantseva and Alain Valette, for their valuable
comments on previous versions of this article. Finally, we would
like to thank Pierre de La Harpe, the referee and the editors for
their careful reading and their interesting suggestions.


\section{Definitions and preliminaries} \label{DefPrelim}

\begin{NB*}
We denote by $\N$ the set $\{0,1,2,\dots\}$ of non-negative integers
and by $\Nast$ (respectively $\Zast$) the sets of positive
(respectively non-zero integers). If $A$ is any ring, then
$A^{\times}$ is the set of invertible elements of $A$. For instance,
$\Z^{\times}=\{-1,1\}$ whereas $\Zast=\Z \setminus \{0\}$.
\end{NB*}

We refer to $\Z$ as the ring of \emph{rational integers} to avoid
confusion with the ring of $m$-adic integers defined in the next
section.

\subsection{The ring of $\m$-adic integers}\label{mAdicSct}

Let $\m \in \Z^*$. Recall that the ring of $\m$-adic integers
$\Z_\m$ is the projective limit of the system
\[
\ldots \to \Z/\m^h\Z \to \Z/\m^{h-1}\Z \to \ldots \to \Z/\m^2\Z \to
\Z/\m\Z
\]
in the category of topological rings, where the arrows are the canonical surjective homomorphisms. This
shows that $\Z_\m$ is compact. This topology is compatible with the
ultrametric distance given, for $\xi \neq \eta$, by
\[
 d_\m(\xi,\eta) = |m|^{-\max\{k\in\N \, : \, \xi-\eta\in \m^k\Z_\m\}} \ .
\]
Now we collect to further use some easy facts about $\m$-adic
integers. Detailed proofs were given in the second-named author's Ph.D. thesis \cite[Appendix
C]{Sta05PhD}\footnote{Note that the second part of Statement (a) was
false there.}. Notice that $\Z_m$ is the zero ring if $\vert m
\vert=1$.
\begin{proposition}\label{mAdic}
Let $\m \in \Z$ such that $|m| \geqslant 2$ and let $\m = \pm
p_1^{k_1} \cdots p_\ell^{k_\ell}$ be its decomposition in prime
factors:
\begin{enumerate}
  \item[(a)] One has an isomorphism of topological rings
  $
    \Z_{\m} \cong \Z_{p_1} \oplus \ldots \oplus \Z_{p_\ell}
  $.
  In particular, if $\m$ is not a power of a prime number, the ring $\Z_{\m}$ has zero divisors.
  \item[(b)] The group of invertible elements of $\Z_\m$ is given by
  $$
    \Z_\m^\times = \Z_\m \setminus (p_1 \Z_\m \cup \ldots \cup
    p_\ell
    \Z_\m).
  $$
  \item[(c)] Any ideal of $\Z_\m$ is principal. More precisely, any non-zero ideal of $\Z_\m$ can
  be written as $p_1^{i_1} \cdots p_\ell^{i_\ell} \Z_\m$ with $i_1, \ldots i_\ell \in \N$.
  \item[(d)] For any $i_1, \ldots i_\ell \in \N$, one has $\Z \cap p_1^{i_1} \cdots p_\ell^{i_\ell} \Z_\m =
  p_1^{i_1} \cdots p_\ell^{i_\ell} \Z$.
\end{enumerate}
\end{proposition}

\begin{definition}\label{gcd}
Let $\m\in\Z$ such that $|m| \geqslant 2$ and let $p_1, \ldots,
p_\ell$ be its prime factors. If $E$ is a subset of $\Z_\m$
containing a non-zero integer, the \emph{greatest common divisor
(gcd)} of the elements of $E$ is the (unique) number $p_1^{i_1}
\cdots p_\ell^{i_\ell}$ (with $i_1, \ldots, i_\ell \in \N$) such
that the ideal generated by $E$ is $p_1^{i_1} \cdots p_\ell^{i_\ell}
\Z_\m$.

If $|m| = 1$, we set by convention $\text{gcd}(\Z_m) = 1$.
\end{definition}

\begin{lemma}\label{Ideals}
Let $\m\in \Z^*$ and let $\m'$ be a divisor of $\m$. Let us write
$\m' = \pm p_1^{j_1} \cdots p_{\ell}^{j_{\ell}}$ and $\m =\pm
p_1^{k_1} \cdots p_\ell^{k_\ell}$ their decomposition in prime
factors ($j_s \leqslant k_s$ for all $s=1, \ldots, \ell$). Let $\pi:
\Z_\m \to \Z_{\m'}$ the ring morphism induced by projections
$\Z/\m^h\Z \to \Z/ (m')^h\Z$ (for $h\geqslant 1$). Then the
following holds:
\begin{enumerate}
    \item[(a)] One has $\pi(\n) = \n$ for any integer $\n$.
    \item[(b)] For any $d=\pm p_1^{i_1} \cdots p_{\ell}^{i_{\ell}}$ with $i_1, \ldots, i_{\ell} \in \N$,
    one has $\pi^{-1}(d\Z_{\m'}) = d \Z_\m$.
\end{enumerate}
\end{lemma}
The ideal $d \Z_m$ is both open and closed in $\Z_m$. More
generally:
\begin{lemma}\label{LemOpenClosedIdeals}
Let $\m\in\Z^*$ and let $m'$ be an integer whose prime divisors
divide $\m$.

There exists $h$ which depends only on $m'$ and such that for any
$\xi,\eta$ in $\Z_m$, the inequality $d_\m(\xi,\eta)<\vert
m\vert^{-h}$ implies $gcd(\xi,m')=gcd(\eta,m')$. In particular, the
set of $m$-adic integers $\xi$ such that $gcd(\xi,m')=d$ is both
open and closed in $\Z_m$.
\end{lemma}

\begin{proof}
Consider $h$ such that $m'$ divides $m^h$. If $d_\m(\xi,\eta)<\vert
m\vert^{-h}$, one has $\xi=k+m^h \mu$ and $\eta=k+m^h\nu$ with $k\in
\Z$ and $\mu,\nu \in \Z_m$ by Proposition \ref{mAdic} (d). Hence
$gcd(\xi,m')=gcd(k,m')=gcd(\eta,m')$.
\end{proof}

\subsection{Marked groups and their topology}\label{MarkGps}

Introductory expositions of these topics can be found in
\cite{Cham00} or \cite{ChGu05}. We recall only the basics needed in this article.

The free group on $k$ generators will be denoted by $\mathbb{F}_k$,
or $\Free_S$ with $S = (s_1, \ldots, s_k)$, if we want to precise
the names of (canonical) generating elements. A \emph{marked group
on $k$ generators} is a pair $(G,S)$ where $G$ is a group and $S =
(s_1, \ldots, s_k)\in G^k$ is a family which generates $G$. 
An \emph{isomorphism of marked groups} on $k$ generators is a group isomorphism which respect the markings. A marked
group $(G,S)$ is endowed with a canonical surjective homomorphism $\phi:
\mathbb{F}_S \to G$, which induces an isomorphism of marked groups
between $\mathbb{F}_S/ \ker \phi$ and $G$. Hence a class of marked
groups for the relation of marked group isomorphism, is represented by a unique quotient of $\mathbb{F}_S$. In
particular if a group is given by a presentation, this defines a
marking on it. The non-trivial elements of ${\cal R} := \ker \phi$
are called \emph{relations} of $(G,S)$. Given $w \in \mathbb{F}_k$
we will often write "$w = 1$ in $G$" or "$w \underset{G}{=} 1$" to
say that the image of $w$ in $G$ is trivial.

Let $w = x_1^{\varepsilon_1} \cdots x_n^{\varepsilon_n}$ be a
reduced word in $\mathbb{F}_S$ (with $x_i \in S$ and $\varepsilon_i
\in \{ \pm 1 \}$). The integer $n$ is called the \emph{length} of
$w$ and denoted $|w|$. If $(G,S)$ is a marked group on $k$
generators and $g \in G$, the \emph{length} of $g$ is
\begin{eqnarray*}
|g|_G & := & \min\{ n: g = s_1 \cdots s_n \text{ with } s_i \in S \sqcup S^{-1} \} \\
                    &  = & \min\{ |w| : w \in \mathbb{F}_S, \ \phi(w) = g  \} \ .
\end{eqnarray*}

Let ${\cal G}_k$ be the set of marked groups on $k$ generators (up
to marked isomorphism). Let us recall that the topology on ${\cal
G}_k$ comes from the following ultrametric distance: for $(G_1, S_1)
\neq (G_2, S_2) \in {\cal G}_k$ we set $d \big( (G_1, S_1), (G_2,
S_2) \big) := e^{-\lambda}$ where $\lambda$ is the length of a
shortest element of $\mathbb{F}_k$ which vanishes in one group and
not in the other one. But what the reader has to keep in mind is the
following characterization of convergent sequences.

\begin{lemma} \rm \cite[Proposition 1]{Sta06a} \it \label{lmecv}
Let $(G_n)_{n}$ be a sequence of marked groups in ${\cal
G}_k$. The sequence $(G_n)_n$ converges if and only if
for any $w \in {\mathbb F}_k$, we have either $w =1$ in $G_n$ for
$n$ large enough, or $w\neq 1$ in $G_n$ for $n$ large enough.
\end{lemma}

The reader could remark that the latter condition characterizes
exactly Cauchy sequences. Sequences converging to a finitely
presented group enjoy a remarkable property\footnote{This property is already used in \cite[Proof of Theorem
6.2]{Gri84}.}.
\begin{lemma}\rm\cite[Lemma 2.3]{ChGu05}\it \label{presFinQuot}
If a sequence $(G_n)_n$ in ${\cal G}_k$ converges to a
marked group $G \in{\cal G}_k$ which is given by a finite
presentation, then, for $n$ large enough, $G_n$ is a marked quotient
of $G$.
\end{lemma}

We address the reader to the given references for proofs.

\subsection{Notation and conventions}\label{NotConv}

We define a family of limits of Baumslag-Solitar groups in the
following way:
\begin{definition}\label{bBSgroups}
For $\m\in\Z^*$ and $\xi\in\Z_\m$, one defines a marked group on two
generators $\bBS{\m}{\xi}$ by the formula
\[
\bBS{\m}{\xi} = \lim_{n \to \infty} BS(\m,\xi_n)
\]
where $(\xi_n)_n$ is any sequence of rational integers such that
$\xi_n$ tends to $\xi$ in $\Z_\m$ and $|\xi_n|$ tends to infinity as
$n$ goes to infinity.
\end{definition}

Notice that $\bBS{\m}{\xi}$ is well defined for any $\xi\in\Z_\m$ by
Theorem 6 of \cite{Sta06a}. Note also that for any $n\in\Z^*$, one
has $\bBS{\m}{\n} \neq BS(\m,\n)$. Indeed, the word $a b^m a^{-1}
b^{-n}$ represents the identity element in $BS(m,n)$, but not in
$\bBS{m}{n}$.

When considered as marked groups, the free group $\F =
\mathbb{F}(a,b)$, Baumslag-Solitar groups, and groups $\BS$ are all
(unless stated otherwise) marked by the pair $(a,b)$.

Another group which plays an important role in this article is:
\[
\Z \wr \Z = \Z \ltimes_t \Z[t,t^{-1}] \cong \Z \ltimes_s
\bigoplus\limits_\Z \Z
\]
where the generator of the first copy of $\Z$ acts on $\Z[t,t^{-1}]$
by multiplication by $t$ or, equivalently, on $\bigoplus_\Z \Z$ by
shifting the indices. This group is assumed (unless specified
otherwise) to be marked by the generating pair consisting of
elements $(1,0)$ and $(0,t^0)$.

The last groups we introduce here are $\Gamma(\m,\n) = \Z
\ltimes_{\frac{\n}{\m}}
\Z[\frac{\text{gcd}(\m,\n)}{\text{lcm}(\m,\n)}]$ ($\m,\n \in \Z^*$)
where the generator of the first copy of $\Z$ acts on
$\Z[\frac{\text{gcd}(\m,\n)}{\text{lcm}(\m,\n)}]$ by multiplication
by $\frac{\n}{\m}$. This group is assumed (unless specified
otherwise) to be marked by the generating pair consisting of
elements $(1,0)$ and $(0,1)$. The latter elements are the images of
$(1,0)$ and $(0,t^0)$ by the homomorphism $\Z \wr \Z \to
\Gamma(\m,\n)$ given by the evaluation $t=\frac{\n}{\m}$; they are
also the images of the elements $a$ and $b$ of $BS(\m,\n)$ by the
homomorphism defined by $a \mapsto (1,0)$ and $b\mapsto (0,1)$.
Observe that the group $\Z \ltimes
\Z[\frac{\text{gcd}(\m,\n)}{\text{lcm}(\m,\n)}]$ acts affinely on
$\Q$ (or $\R$) by $(1,0) \cdot x = \frac{\n}{\m}x$ and $(0,y)\cdot x
= x+y$.

We introduce the homomorphism $\sigma_a: \F \to \Z$ defined by
$\sigma_a(a) = 1$ and $\sigma_a(b) = 0$. It factors through all
groups $BS(\m,\n)$, $\BS$, $\Z \wr \Z$ and $\Gamma(\m,\n)$. The
induced morphisms are also denoted by $\sigma_a$. To end this
section we define the homomorphism $\bar{} : \F \to \F$ given by
$\bar a = a$ and $\bar b = b^{-1}$. Note that it induces
homomorphisms $\bar{} : BS(\m,\n) \to BS(\m,\n)$ for $\m,\n\in\Z^*$
and $\bar{} : \bBS{\m}{\xi} \to \bBS{\m}{\xi}$ for
$\m\in\Z^*,\xi\in\Z_\m$.

\section{Parametrization of limits} \label{Sctiff}

This section is devoted to the map
\[
 \BBS_m :  \Z_m  \longrightarrow \G ;\, \xi  \longmapsto \BS \ .
\]
First, we characterize its lack of injectivity, thus giving a
classification of groups $\BS$ up to isomorphism of marked groups
(Theorem \ref{ThmBSm}).
 Second, we prove that it is continuous (Corollary \ref{Embedding}).

\subsection{Lack of injectivity}
Let us summarize the key points of this set-theoretic part. Our aim
is the classification (up to isomorphism of marked groups) of the
groups $\BS$ marked by the pair $(a,b)$:
\begin{theorem} \label{ThmBSm}
Let $\m \in\Z^{\ast}$ and $\xi,\eta\in\Z_{\m}$. The equality of
marked groups $$\overline{BS}(\m,\xi)=\overline{BS}(\m,\eta)$$ holds
if and only if there is some $d\in\N^{\ast}$ such that
$gcd(\xi,\m)=gcd(\eta,\m)=d$ and the images of $\xi/d$ and $\eta/d$
in $\Z_{\m/d}$ are equal.
\end{theorem}

As a consequence, we obtain the injectivity of $\BBS_m$, when
restricted  to the set invertible $m$-adic integers.
\begin{corollary}\label{CorBSm}
Let $\m\in\Z^*$ and let $\xi, \eta\in\Z_\m$ with $\xi \neq \eta$. If
no prime factor of $\m$ divides both $\xi$ and $\eta$, then one has
$\bBS{\m}{\xi} \neq \bBS{\m}{\eta}$. Hence the map $\BBS_m$ becomes
injective when restricted to invertible $m$-adic integers.
\end{corollary}

\begin{sproof}{Proof of Corollary \ref{CorBSm}}
By Theorem \ref{ThmBSm}, we may suppose that
$\text{gcd}(\m,\xi)=\text{gcd}(\m,\eta) \linebreak =d.$ Then $d=1$
by assumption and $\pi(\xi)=\xi \neq \eta=\pi(\eta)$ where $\pi:
\Z_m \twoheadrightarrow \Z_{m/d}$ is the canonical ring morphism. By
Theorem \ref{ThmBSm}, one gets $\overline{BS}(\m,\xi) \neq
\overline{BS}(\m,\eta)$.
\end{sproof}

To prove Theorem \ref{ThmBSm}, we first need to characterize the converging
sequences of Baumslag-Solitar groups marked by the standard pair
$(a,b)$. This is done with the following statement:
\begin{theorem} \label{thmiffconv}
Let $\m\in\Z^{\ast}$ and let $(\xi_{\n})_{\n}$ be sequence of
rational integers such that $\vert \xi_n \vert$ tends to infinity as
$n$ goes to infinity.
 The sequence $(BS(\m,\xi_{\n}))_{\n}$ converges in $\mathcal{G}_2$ if and only if the
 two following conditions hold:
 \begin{itemize}
 \item[$(i)$] there is $d\in \Nast$ such that $gcd(\m,\xi_{\n})=d$ for all $\n$ large enough;
 \item[$(ii)$] $(\xi_n/d)_n$ converges in $\Z_{m/d}.$
 \end{itemize}
\end{theorem}
One implication of Theorem \ref{thmiffconv} was proved by the
second-named author \cite{Sta06a} for $d=1$. The purpose of the
following proposition is to generalize it to any $d$:
\begin{proposition} \label{propsuff}
Let $\m\in\Z^{\ast}$ and let $(\xi_n)_{n}$ be a sequence of rational
integers such that $\vert \xi_n \vert$ tends to infinity as $n$ goes
to infinity. If $(\xi_n)_n$ defines a converging sequence in
$\Z_{\m}$ then the sequence of marked groups $\left( BS(md,\xi_n d)
\right)_{n}$ converges in $\mathcal{G}_2$ for any $d \in \Zast$.
\end{proposition}

We obtain this proposition by performing minor technical
modifications on \cite[Theorem 6]{Sta06a} and the related lemmas.
Nevertheless, the process being non completely obvious, we give a a
detailed proof of it. As the proof of Proposition \ref{propsuff} is
more complex than the first two statements, we postponed it to the
end of this subsection.

We first show that Theorem \ref{thmiffconv} implies Theorem
\ref{ThmBSm}.

\begin{sproof}{Proof of Theorem \ref{ThmBSm}}
 Choose a sequence of rational integers $(\xi_{\n})_{\n}$ such that
 $$
  \xi_{2 \n} \underset{\Z_{\m}}{\longrightarrow} \xi,\,\xi_{2 \n+1} \underset{\Z_{\m}}{\longrightarrow} \eta
  \mbox{ and } \vert \xi_{\n} \vert \longrightarrow \infty \mbox{ as } \n \mbox{ goes to infinity}.
 $$
 One has $\overline{BS}(\m,\xi) = \overline{BS}(\m,\eta)$ if and only if the sequence $(BS(\m,\xi_{\n}))_n$ converges.
  By Theorem \ref{thmiffconv} it is equivalent to have $\text{gcd}(\m ,\xi_{\n})=d$ for $\n$ large enough (for some $d$ in $\N^{\ast}$)
   and the sequence $(\pi(\xi_n/d))_n$ converging in $\Z_{m/d}$ (where $\pi: \Z_m \twoheadrightarrow \Z_{m/d}$ is the canonical ring morphism).
  Finally, the first condition is equivalent to $\text{gcd}(\m,\xi)=d=\text{gcd}(\m,\eta)$ by Lemma
  \ref{LemOpenClosedIdeals} and the second condition boils down to $\pi(\xi/d)=\pi(\eta/d)$.
 \end{sproof}

We now prove that Proposition \ref{propsuff} and the following lemma
imply Theorem \ref{thmiffconv}. This lemma is actually independent
and is only used for the converse implication.

\begin{lemma} \label{lemneqd} Let $m_i,d_i,k_i\in\Z^{\ast}$ for
$i=1,2$ such that
$$m_1d_1=m_2 d_2,\,\vert k_2d_2 \vert \neq 1,\, gcd(m_2,k_2)=1 \mbox{ and }
d_1 \mbox{ does not divide } d_2.
$$

Then, the distance between $BS(m_1d_1,k_1d_1)$ and $BS(m_2d_2,k_2
d_2)$ in $\mathcal{G}_2$ is not less than $e^{-\delta}$ with
$\delta=10+2d_1 m_1^2.$
\end{lemma}

\begin{sproof}{Proof of Lemma \ref{lemneqd}}
Consider $r=a^2b^{d_1 m_1^2}a^{-2}b$ and let $w=r \overline{r}$. On
one hand, we have $r=b^{d_1 k_1^2+1}$ in $BS(m_1d_1,k_1d_1)$, which
implies $w=1$ in $BS(m_1d_1,k_1d_1).$ On the other hand,
$r=ab^{m_1d_2k_2}a^{-1}b$ in $BS(m_2 d_2,k_2 d_2)$. As $m_2$ and
$k_2$ are coprime integers, we deduce that $m_2 d_2$ divides $m_1
d_2 k_2$ if and only if $d_1$ divides $d_2$. Under the assumptions
of the lemma, the writing
$ab^{m_1d_2k_2}a^{-1}bab^{-m_1d_2k_2}a^{-1}b^{-1}$ is then a reduced
form for $w$ in $BS(m_2d_2,k_2 d_2)$. By Britton's Lemma, we have $w
\neq 1$ in $BS(m_2d_2,k_2 d_2)$. As $\vert w \vert=10+2d_1 m_1^2$,
we get the conclusion.
\end{sproof}

\begin{sproof}{Proof of Theorem \ref{thmiffconv}}
Proposition \ref{propsuff} implies immediately that conditions (i)
and (ii) are sufficient.

Let us show that conditions (i) and (ii) are necessary. To do so, we
assume that the sequence $(BS(\m,\xi_{\n}))_{\n}$ converges in
$\mathcal{G}_2$. If (i) does not hold, we can find two subsequences
$(\xi'_{\n})_{\n}$ and $(\xi''_{\n})_{\n}$ of $(\xi_{\n})_{\n}$ such
that $gcd(m,\xi'_{\n})=d_1, gcd(m,\xi''_{\n})=d_2, \vert \xi''_{\n}
\vert>1$ for all $\n$ and $d_1$ does not divide $d_2$. Then Lemma
\ref{lemneqd} clearly shows that $(BS(\m,\xi_{\n}))_{\n}$ is not a
converging sequence in $\mathcal{G}_2$, a contradiction. Hence (i)
holds; let $d$ be a rational integer satisfying it. The marked
subgroup $\Gamma_{\xi_{\n},d}$ of $BS(m,\xi_n)$ generated by
$(a,b^d)$ is equal to $BS(m/d,\xi_n/d)$ endowed with its standard
marking $(a,b)$. The sequence of the $BS(m/d,\xi_n/d)$'s is then
also converging in $\mathcal{G}_2$. By Theorem 3 of \cite{Sta06a},
the sequence of rational integers $(\xi_{\n}/d)_{\n}$ defines a
converging sequence in $\Z_{m/d}$. Hence (ii) holds.
\end{sproof}

We now turn to the proof of Proposition \ref{propsuff}, which occupies
the end of this Subsection. Because of
Lemma \ref{lmecv}, it boils down to show that for any $w$ in $\F$ we
have the implication: if $w$ reduces to $1$ in $G_n=BS(md,\xi_nd)$
for infinitely many $n$ then $w$ reduces to $1$ in $G_n$ for all $n$
large enough. Our strategy is to prove that under the hypothesis of
Proposition \ref{propsuff}, any word $w$ undergoes the same sequence
of cancelations in $G_n$ for all $n$ large enough.

Let $(G,(a,b))$ be a marked group on two generators. We call
\emph{$a$-length} of an element $g\in G$ the minimal number of
letters $a,a^{-1}$ occurring in a word on $\{a,a^{-1},b,b^{-1}\}$
which represents $g$.

The precise statement is:

\begin{lemma} \label{LemSimCanc}
Let $w=b^{e_0} a^{\varepsilon_1} b^{e_1} \cdots a^{\varepsilon_h}
b^{e_h}$ with $\varepsilon_j =\pm 1$ for $j=1,\dots,h$. Let $0 \le t
\le h/2$ and let $C$ be a class modulo $m^t$. Assume that there are
infinitely many $n\in C$ such that the $a$-length of (the image of)
$w$ in $BS(md,nd)$ is at most $h-2t$. Then there exist
$\delta_1,\ldots,\delta_{h-2t} \in \{\pm 1\}$ and polynomial
functions $\alpha_0, \ldots, \alpha_{h-2t}$, depending only on $w$
and $C$, such that
\[
w \underset{BS(md,nd)}= \bar w(n) := b^{\alpha_0(n)} a^{\delta_1}
b^{\alpha_1(n)} \cdots a^{\delta_{h-2t}} b^{\alpha_{h-2t}(n)}
\]
for all $n\in C$ with $|n|$ large enough.
\end{lemma}

Proposition \ref{propsuff} then follows:

\begin{sproof}{Proof of Proposition \ref{propsuff}}
Assume that $w=b^{e_0} a^{\varepsilon_1} b^{e_1} \cdots
a^{\varepsilon_h} b^{e_h}$ is trivial in $BS(md,\xi_n d)$ for
infinitely many $n$. The sum $\varepsilon_1+\dots+\varepsilon_h$ has
to be zero. Thus $h$ is even and we can set $t=h/2$. By hypothesis,
there exists a class modulo $m^t$, say $C$, such that $\xi_n \in C$
for all $n$ large enough. We apply Lemma \ref{LemSimCanc}: there
exists a polynomial function $\alpha=\alpha_0$ depending only on $w$
and $C$ such that $w$ reduces to $\bar w(n)=b^{\alpha(\xi_n)}$ in
$BS(md,\xi_nd)$ for all $n$ large enough. As $\alpha(\xi_n)$ is zero
for infinitely many $n$, the polynomial function $\alpha$ is the
zero function. Hence $w=1$ in $BS(md,\xi_nd)$ for all $n$ large
enough, which proves that $(BS(md,\xi_nd))_n$ converges in $\G$ by
Lemma \ref{lmecv}.
\end{sproof}

It only remains to prove Lemma \ref{LemSimCanc}. Reductions of $w$
in $BS(md,\xi_nd)$ are ruled by the congruence classes of its $b$
exponents modulo $md$ and $\xi_nd$. By Britton's Lemma, we have
indeed: $w=b^{e_0} a^{\varepsilon_1} b^{e_1} \cdots
a^{\varepsilon_h} b^{e_h}$ is not reduced in $BS(md,\xi_nd)$ if and
only if there is some $i$ in $\{1,\dots,h\}$ such that
\begin{itemize} \item[$(1)$] either $\varepsilon_i=-\varepsilon_{i+1}=1$ and $e_i
\equiv 0 \, (\text{mod }md)$ \item[$(2)$] or
$\varepsilon_i=-\varepsilon_{i+1}=-1$ and $e_i \equiv 0 \,(\text{mod
}\xi_nd)$.
\end{itemize}
In both cases, the word $w$ reduces in $BS(md,\xi_nd)$ to a word
$w'$ such that its length with respect to $\{a^{-1},a\}$ decreases
by $2$. We define below integers
$r_0(\xi_n),\dots,r_t(\xi_n),s_0(\xi_n),\dots,s_t(\xi_n)$ to keep
track of the $b$ exponents in the reductions $w',w'',\dots,w^{(t)}$
of $w$ in $BS(md,\xi_nd)$ along a sequence of $t$ cancelations of
type $(1)$ or $(2)$. These integers will be crucially used in the
proof of Lemma \ref{LemSimCanc} to give a description of such $b$
exponents as functions of $n$. The integers $r_i(\xi_n)$
(respectively $s_i(\xi_n)$ modulo $m$) are shown to depend only on
the congruence class of $\xi_n$ modulo $m^t$.

\begin{lemma}\label{Prelim1}\rm\cite[Lemma 4]{Sta06a} \it
    Fix $\m\in\Z^{\ast}$. We define recursively two sequences $s_0, s_1,\ldots$ and $r_0,r_1,\ldots$ of functions from $\Za$ to $\Z$ by
    \begin{enumerate}
        \item[(i)]  $r_0(n)=0$ and $s_0(n) = 1$ for all $n$;
        \item[(ii)] $s_{i-1}(n) \n = s_i(n) \m + r_i(n) $ and $0 \leqslant r_i(n) < \m$ for $i \geqslant 1$ (Euclidean division).
    \end{enumerate}
    Then, for any $n,n'$ in $\Z^{\ast}$ and $t \ge 1$ such that $n' \equiv n \, (\text{mod } \m^{t})$, we have
     $r_i(n)=r_i(n')$ and $s_i(n) \equiv s_i(n') \ (\text{mod } \m^{t-i})$ for all $0 \leqslant i \leqslant t$.
\end{lemma}

\begin{remark}\label{remPol} Let us observe that the functions $s_0,\dots,s_t$
are polynomial functions when restricted to a given class $C$ modulo
$m^t$. Let $c \in \{0,\dots,m^t-1\}$ and let $C$ be the set of
rational integers $n$ such that $n \equiv c \,(\text{mod } m^t)$.
Define recursively the polynomial $P_{i,C}(X)$ by $P_{0,C}(X)=1$ and
$P_{i,C}(X)=XP_{i-1,C}(X)-\frac{r_i(n)}{m}$ with $n \in C$ and $i
\ge 1$. Clearly, $mP_{i,C}(X)=mX^{i}-r_1(n)X^{i-1}-\dots-r_i(n)$.
The previous lemma implies that $P_{i,C}$ is a polynomial whose
coefficients do not depend on $n \in C$. As
$s_i(n)=P_{i,C}(\frac{n}{m})$ for all $0 \le i \le t$ and all $n\in
C$, $s_i$ is a polynomial function of degree $i$ on $C$.

\end{remark}

The following lemma describes how the $b$ exponents in $w$ transform
through a cancelation of type $(1)$ or $(2)$ in $BS(md,nd)$. These
exponents are put into a particular form by means of the $s_i$ and
this form is shown to be (fortunately) preserved under reductions.
The result coincides with \cite[Lemma 5]{Sta06a} when $d=1$. The
proof is very similar: it is provided for completeness.

\begin{lemma} \label{lemard}
Fix $m,d\in\Z^{\ast}$ and $t \ge 1$. Let $C$ be a class modulo
$m^t$. Let $k_0,\dots,k_t \in \Z$ and let $\alpha: C \rightarrow \Z$
be the function defined by:
$$
\alpha(n)=k_0+k_1 \n d+k_2 s_1(n) \n d+\dots+k_t s_{t-1}(n) \n d
$$

where $s_0, s_1, \ldots$ are given by Lemma \ref{Prelim1}. Let us
also take $r_0, r_1, \ldots$ as in Lemma \ref{Prelim1}.

 \begin{itemize}
 \item[$(i)$] The class of $\alpha(n)$ modulo $md$ does not depend
 on $n \in C$.
 If $\alpha(n) \equiv 0 \,(\text{mod }md)$ for some $n \in C$, then $d$ divides $k_0$ and
 we get $ab^{\alpha(n)}a^{-1} =b^{\beta(n)}$ in $BS(md, nd)$ for all $n\in C$,
  with
  $$
  \beta(n)= l_1 n d+l_2 s_1(n) n d+\dots+l_{t+1} s_t(n) n d
  $$
  and
  $$
  l_1=\frac{1}{\m}\left(\frac{k_0}d+k_1r_1(n)+\dots+k_t r_t(n)\right),\,l_i=k_{i-1} \mbox{ for } 2 \le i \le t+1.
  $$
  \item[$(ii)$] We have
 either $\alpha(n) \equiv 0$ $(\text{mod } \n d)$ for all $n \in
C$ such that $|n|>k_0$
 or $\alpha(n) \not \equiv 0$ $(\text{mod } \n d)$ for all $n \in
C$ such that $|n|>k_0$.

 In the first case
 we get $a^{-1}b^{\alpha(n)}a =b^{\beta(n)}$ in $BS(\m d, \n d)$
 for all $n \in C$ such that $|n|>k_0$,  with
\[
  \beta(n)= l_0  d+l_1 \n d+l_2 s_1(n) \n d+\dots+l_{t-1} s_{t-2}(n) \n d
\]
  and
 \[
  l_0=k_1 \m-k_2 r_1(n)-\dots-k_{t} r_{t-1}(n),\,l_i=k_{i+1} \mbox{ for } 1 \le i \le t-1 \ .
 \]
 \end{itemize}
Moreover, the $l_i$'s are constant rational integers in both cases.
\end{lemma}

\begin{proof}
$(i)$ Lemma \ref{Prelim1} ensures that the class of $s_i(n)$ modulo
$m$ does not depend on $n \in C$ for all $i=1,\dots,t-1$. Hence the
class of $\alpha(n)$ modulo $md$ does not depend on $n \in C$.
Assume now that $ \alpha(n) \equiv 0 \, (\text{mod } md) $ for some
$n \in C$. Then, we have $k_0\equiv 0 \, (\text{mod } \m d)$,
$ab^{\alpha(n)}a^{-1} =b^{\n\alpha(n)/\m}$ in $BS(\m d, \n d)$ with
\[
  \alpha(n) = \left(\frac{k_0}d + k_1 r_1(\n) + \cdots + k_t r_t(\n) + k_1 s_1(\n)\m  + \cdots + k_t s_t(\n) \m \right)d
\]
since $s_{i-1}(n)n = s_{i}(n) m + r_{i}(n)$ for all $i\ge 1$. Hence, we
obtain
\begin{eqnarray*}
 \frac{\n\alpha(n)}{\m} & = & \frac 1\m \left(\frac{k_0}d + k_1 r_1(\n) + \cdots + k_t r_t(\n) \right) nd + k_1 s_1(\n)\n d  + \cdots + k_t s_t(\n) \n d  \\
 & = & l_1 nd + l_2 s_1(\n)\n d  + \cdots + l_{t+1} s_t(\n) \n d = \beta(n) \ .
\end{eqnarray*}
By Lemma \ref{Prelim1}, $r_i(n)$ does not depend on $n \in C$ for
all $i = 1,\ldots, t$. Consequently, the rational integers
$l_1,\dots,l_{t+1}$ do not depend on $n \in C$.

\medskip

$(ii)$ Assume that $\vert n \vert> \vert k_0 \vert$. We have
$\alpha(n) \equiv 0 \, (\text{mod } nd)$ if and only if $ k_0=0$.
Suppose now that it is the case. We have then $ab^{\alpha(n)}a^{-1}
=b^{\m\alpha(n)/\n}$ in $BS(\m d, \n d)$ with
\begin{eqnarray*}
  \frac{\m\alpha(n)}{\n} & = & k_1 \m d + k_2 s_1(\n)\m d + \cdots +
  k_t s_{t-1}(\n)\m d \\
  & = & k_1 \m d - k_2 r_1(\n) d - \cdots - k_t r_{t-1}(\n) d + k_2 s_0(\n)\n d
   + \cdots + k_t s_{t-2}(\n) \n d \\
  & = & l_0  d+l_1 \n d+l_2 s_1(n) \n d+\dots+l_{t-1} s_{t-2}(n) \n d = \beta(n)
\end{eqnarray*}
since $s_{i}(n) m = s_{i-1}(n)n - r_{i}(n)$ for all $i\ge 1$. By
Lemma \ref{Prelim1} again, the rational integers $l_0,\dots,l_{t-1}$ do
not depend on $n \in C$.
\end{proof}

\medskip

We are eventually able to prove Lemma \ref{LemSimCanc}.

\begin{sproof}{Proof of Lemma \ref{LemSimCanc}}
Set $G_n=BS(md,nd)$. We show by induction on $t\ge 0$ that, provided
the assumption in the Lemma is satisfied, the functions $\alpha_i: C
\longrightarrow \Z$ exist and satisfy:
\begin{equation} \label{EqPolySi}\tag{\textasteriskcentered}
\left\{
\begin{array}{rcll}
 \alpha_i(n) & = & k_{0,i} & \text{ if } t=0   \\
 \alpha_i(n) & = &
   k_{0,i}+k_{1,i}nd+k_{2,i}s_2(n)nd+\dots+k_{t,i}s_{t-1}(n)nd &
\text{ if } t>0
\end{array}
\right. \ ,
\end{equation}
where the $s_j$ are the functions defined in Lemma \ref{Prelim1} and
$k_{0,i},\dots,k_{t,i}$ are integers that depend only on $w$ and
$C$. Functions of the form (\ref{EqPolySi}) depend only on $w$ and
$C$. These functions are polynomial functions of $n \in C$ by Remark
\ref{remPol}.

\textbf{Case $t=0$:} It suffices to take $\alpha_i(n) = e_i$ to
satisfy (\ref{EqPolySi}) and $\bar w(n)=w$ for all $n$.

\textbf{Induction step ($0<t\le h/2$):} We denote by $C'$ the class
modulo $m^{t-1}$ defined by $C$. By assumption, there exists an
infinite subset $I\subseteq C \subset C'$ such that the $a$-length
of $w$ in $BS(md,nd)$ is at most $h-2t$ for all $n \in I$. By
induction hypothesis, there exist $\delta'_1, \ldots,
\delta'_{h-2t+2} \in \{\pm 1\}$ and polynomial functions $\alpha'_0,
\ldots, \alpha'_{h-2t + 2}$ satisfying
\begin{equation*}
\left\{
\begin{array}{rcll}
 \alpha'_i(n) & = & k_{0,i} & \text{ if } t=1   \\
 \alpha'_i(n) & = &
   k_{0,i}+k_{1,i}nd+k_{2,i}s_2(n)nd+\dots+k_{t-1,i}s_{t-2}(n)nd &
\text{ if } t>1
\end{array}
\right. \ ,
\end{equation*}
and such that $w = w'(n) := b^{\alpha'_0(n)} a^{\delta'_1}
b^{\alpha'_1(n)} \cdots a^{\delta'_{h-2t+2}}
b^{\alpha'_{h-2t+2}(n)}$ in $BS(md,nd)$ for all $n\in C'$ with $|n|$
large enough.

Hence, for all $n\in I$ with $|n|$ large enough, the word $w'(n)$ is
not reduced in $BS(md,nd)$, since its $a$-length in this group is at
most $h-2t$. For such $n$'s, Britton's Lemma gives an index
$i=i(n)\in \{1,\ldots, h-2t+1\}$ such that:
\begin{enumerate}
 \item[(a)] either $a^{\delta'_i}b^{\alpha'_i(n)}a^{\delta'_{i+1}} =
 ab^{\alpha'_i(n)}a^{-1}$ and $\alpha'_i(n) \equiv 0$ (mod $md$);
 \item[(b)] or $a^{\delta'_i}b^{\alpha'_i(n)}a^{\delta'_{i+1}} =
 a^{-1}b^{\alpha'_i(n)}a$ and $\alpha'_i(n) \equiv 0$ (mod $nd$).
\end{enumerate}
Therefore, there exists an infinite subset $I'\subseteq I\subseteq
C$ and a \textbf{fixed} index $i$, such that one of the cases
(a),(b) occurs for all $n\in I'$. Applying Lemma \ref{lemard}, we
see that a reduction will occur at index $i$ for all $n\in C$ with
$|n|$ large enough. Hence, we get $a^{\delta'_i}b^{\alpha'_i(n)}
a^{\delta'_{i+1}} = b^{\beta(n)}$ in $BS(md,nd)$, where
\[
 \beta(n) = l_0 d + l_1 n d+l_2 s_1(n) n d+\dots+l_{t+1} s_t(n) n d \ ,
\]
for all $n\in C$ with $|n|$ large enough (two last terms vanishing
if we were in case (b)). Finally, we are done by setting
\[
\begin{array}{c}
\alpha_j(n):=\left\{ \begin{array}{ccc} \alpha'_j(n) &\text{ if }&
j<i-1;\\
\alpha'_{i-1}(n)+\beta(n)+\alpha'_{i+1}(n)&\text{ if }& j=i-1;\\
\alpha'_{j-2}(n) &\text{ if }& j \ge i;
\end{array}\right.
\end{array}
\]
for $j=0, \ldots, h-2t$.
\end{sproof}

\subsection{Continuity} \label{SctInj}

We now turn to the topologist's point of view. Namely, we are going
to prove the following analogue of Theorem \ref{thmiffconv}:

\begin{theorem} \label{thmiffconvZm}
Let $\m$ be in $\Z^{\ast}$ and let $(\xi_{\n})_{\n}$ be a sequence
in $\Z_m$.
 The sequence $(\overline{BS}(\m,\xi_{\n}))_{\n}$ converges in $\mathcal{G}_2$ if and only if the following conditions both
 hold:
 \begin{itemize}
 \item[$(i)$] there is $d\in\Nast$ such that $gcd(\m,\xi_{\n})=d$ for all $\n$ large enough;
 \item[$(ii)$] the images of $\xi_n/d$ define a converging sequence in $\Z_{m/d}.$
 \end{itemize}
 In particular, if $(\xi_n)_n$ converges to $\xi$ in $\Z_m$, then the sequence
 $(\overline{BS}(m,\xi_n))_n$ converges to $\BS$ in $\G$.
\end{theorem}

The following corollary is immediate:

\begin{corollary}\label{Embedding}
For all $\m \in \Z^*$, the map $\overline{BS}_\m : \Z_\m \to \G$ is
continuous.
\end{corollary}

For $|\m| \geqslant 2$, note that if we endow $\Z$ with the
$\m$-adic ultrametric, the analogue map
$$
\Z^* \to \G  ; \ \n \mapsto BS(\m,\n)
$$
is nowhere continuous. Indeed, one has $\n + m^k$ tends to $n$ as
$k$ tends to infinity, while $BS(\m,\n+m^k) \to \bBS{\m}{\n}$. We
recall that one has $\bBS{\m}{\n} \neq BS(\m,\n)$, since the word
$ab^\m a^{-1}b^{-\n}$ defines the trivial element in $BS(\m,\n)$ but
not in $\bBS{\m}{\n}$.

\begin{sproof}{Proof of Theorem \ref{thmiffconvZm}}
Let $(\xi_n)_n$ be a sequence of $m$-adic integers. By an easy
argument of diagonal extraction, we can find a rational integer
sequence $(\eta_n)_n$ such that:
\begin{itemize}
\item[$(1)$] $d(BS(m,\eta_n),\overline{BS}(m,\xi_n))
< \frac{1}{n}$;
\item[$(2)$] $d_m(\eta_n,\xi_n)< \frac{1}{n}$;
\item[$(3)$]$\vert \eta_n \vert \ge n$.
\end{itemize}
The sequence $(\overline{BS}(m,\xi_n))_n$ converges if and only if
$(BS(m,\eta_n))_n$ converges because of inequality $(1)$. By Theorem
\ref{thmiffconv}, the sequence $(BS(m,\eta_n))_n$ converges if and
only if there is some $d$ in $\Z^*$ such that $gcd(\eta_n,m)=d$ for
all $n$ large enough and $(\eta_n/d)_n$ converges in $\Z_{m/d}$.
 By $(2)$ and Lemma \ref{LemOpenClosedIdeals}, we have $gcd(\xi_n,m)=d$ for all $n$ large enough if and only if $gcd(\eta_n,m)=d$ for all $n$ large
enough. Since $d_{m/d}(\eta_n/d,\pi(\xi_n/d)) \le
d_m(\eta_n,\xi_n)$, we deduce from $(2)$ that the sequence
$(\eta_n/d)_n$ converges in $\Z_{m/d}$ if and only if
$(\pi(\xi_n/d))_n$ converges in $\Z_{m/d}$.

Now, let us turn to the particular case $\xi_n \to \xi$ in $\Z_\m$.
We get condition (i) by Lemma \ref{LemOpenClosedIdeals} and
condition (ii) is obvious. The sequence $(\overline{BS}(\xi_n,m))_n$
is thus converging. Moreover, its limit is $\bBS{\m}{\xi}$, since it
does not change if one intercalates a $\xi$-constant subsequence in
$(\xi_n)_n$.
\end{sproof}

\smallskip

We now ``particularize to the case of invertible elements'' and show
that, in this case, the $\overline{BS}$ groups form the set of
accumulation points of the $BS$ groups. Precise statements are as
follows:

\begin{definition}
For $\m\in\Z^*$, we define:
\begin{eqnarray*}
  X_\m & = & \big\{ BS(\m,\n) \ : \ \n \text{ is relatively prime to } \m \big\} \ ; \\
  Z_\m^\times & = & \big\{ \bBS{\m}{\xi} \ : \ \xi \in \Z_\m^\times \big\} \ .
\end{eqnarray*}
By convention, we say that $Z_{\pm 1}^\times$ is empty.
\end{definition}

\begin{corollary}\label{Bdy}
For all $\m\in\Z^*$, the boundary of $X_\m$ in ${\cal G}_2$ is
$Z_\m^\times$. It is homeomorphic to the set of invertible $\m$-adic
integers.
\end{corollary}

\begin{proof} Theorem 3 of \cite{Sta06a} implies that the elements of $\overline{X_\m}$ are the $BS(\m,\n)$'s
with $\n$ relatively prime to $\m$ and the $\bBS{\m}{\xi}$ with $\xi
\in \Z_\m^\times$. One sees easily that the $BS(\m,\n)$'s are
isolated points in $\overline{X_\m}$ (consider the word $a b^\m
a^{-1} b^{-\n}$). The equality $\partial X_\m = Z_\m^\times$ follows
immediately. The second statement is a direct consequence of
Corollary \ref{Embedding}.
\end{proof}

\section{Actions and structure of the limits}\label{Actions}

It is well-known that Baumslag-Solitar groups act on their Bass-Serre trees
and (affinely) on the real line or on $\Q$. In this section, we study such
actions (which are ``trivial'' in the affine case) and the structure of the
limits $\bBS{\m}{\xi}$.

\subsection{A common quotient}

We recall that $\Gamma(\m,\n) = \Z \ltimes_\frac{\n}{\m}
\Z[\frac{\text{gcd}(\m,\n)}{\text{lcm}(\m,\n)}]$ acts affinely on
$\R$ and that it is a marked quotient of both $BS(\m,\n)$ and $\Z
\wr \Z$ (see Section \ref{NotConv}). Notice that these groups are distinct
from the groups $\Gamma_{\xi_n,d}$ introduced in Section \ref{Sctiff}.

The following lemma enables us to define converging sequences of
homomorphisms induced by the same endomomorphism of $\Fm$.
\begin{lemma} \label{LemLimHom}
Let $(G_n)_n$ and $(H_n)_n$ be sequences of marked groups and let
$G$ and $H$ be their limits in $\Gm$. Let $\phi$ be an endomorphism of
$\Fm$ which induces an homomorphism from $G_n$ to $H_n$ for all $n$.
Then $\phi$ induces an homomorphism from $G$ to $H$. Moreover, if
$\phi : G_n \longrightarrow H_n$ is injective for infinitely many $n$,
then $\phi : G \longrightarrow H$ is injective.
\end{lemma}

\begin{proof}
Let $w \in \Fm$ be such that $w \eg{G} 1$. For all $n$ large enough
we have $w \eg{G_n} 1$ (Lemma \ref{lmecv}), and hence $\phi(w)
\eg{H_n} 1$. Consequently, $\phi(w) \eg{H}1$, which shows that
$\phi$ induces an homomorphism from $G$ to $H$. Assume $\phi$ is
injective for infinitely many $n$ and consider $w \in \Fm$ such that
its image in $G$ belongs to $\ker \left (\phi : G \longrightarrow H
\right)$. There are infinitely many $n$ such that the following both
hold :
\begin{itemize}
\item[$(i)$] $\phi : G_n \longrightarrow H_n$ is injective;
\item[$(ii)$] $\phi(w) \eg{H_n} 1$.
\end{itemize} It follows that $w \eg{G_n}1$
holds for infinitely many $n$. As a consequence, $w \eg{G}1$.
\end{proof}

\begin{proposition} \label{propquot}
For any $\m \in \Z^\ast$ and $\xi \in \Z_{\m}$, the morphism of
marked groups \\${q:\F \twoheadrightarrow \Z\wr\Z}$
 factors through a morphism $q_{\m,\xi}: \overline{BS}(\m,\xi) \twoheadrightarrow \Z \wr \Z$.
\end{proposition}

\begin{proof} Take $(\xi_n)_{n}$ a
sequence of rational integers such that $\vert \xi_n \vert$ tends to
infinity and $\xi_n$ tends to $\xi$ in $\Z_{\m}$ as $n$ goes to
infinity. By Definition \ref{bBSgroups}, the sequence
$(BS(\m,\xi_{\n}))_n$ converges to $\overline {BS}(\m,\xi)$. We now
claim that $\Gamma(\m,\xi_{\n})$ tends to $\Z \wr \Z$ as $n$ goes to
infinity, a statement analogous to Theorem 2 in \cite{Sta06a}. As
each $\Gamma(\m,\xi_{\n})$ is a marked quotient of $\Z\wr\Z$, we
only have to prove that any non-trivial element $(\sigma, P(t)) \in
\Z\wr\Z$ defines eventually non-trivial elements in the
$\Gamma(\m,\xi_{\n})$'s. This is obvious if $\sigma \neq 0$. Assume
$\sigma=0$. The Laurent polynomial $P$ has only finitely many roots,
so that $P(\xi_n/m)\neq 0$ for $n$ sufficiently large. This tells us
that $(\sigma,P)$ is eventually non-trivial in the
$\Gamma(\m,\xi_{\n})$'s and proves the claim.

Hence, one has the following ``diagram'':
$$
\xymatrix{
BS(\m,\xi_{\n}) \ar@{->>}[d] \ar[r]^{n \to \infty} & \overline{BS}(\m,\xi) \\
\Gamma(\m,\xi_{\n}) \ar[r]^{n \to \infty} & \Z \wr \Z }$$ As
vertical arrows are induced by the identity of $\F$, the result
follows from the previous lemma.
\end{proof}

\subsection{Affine actions} \label{Aff}

Let $R$ be a ring and $\text{Aff}(R) \cong R^\times\ltimes R$ be the affine
group on $R$. For any $\alpha \in R^\times$, it is straightforward to see that
formulae $a\cdot x = \alpha x$ and $b\cdot x = x+1$ define a homomorphism
$\Z\wr\Z \to \text{Aff}(R)$. Thus Proposition \ref{propquot} immediately
gives the following:
\begin{proposition}\label{AffineAction}
Let $R$ be a ring and let $\alpha\in R^\times$. For any $\m\in \Z^\ast$ and any
$\xi\in \Z_\m$, the group $\BS$ acts affinely on $R$ via formulae:
\[
 a\cdot x = \alpha x \ \text{ and } \ b\cdot x = x+1 \ .
\]
\end{proposition}

\subsection{Actions on trees}\label{Tree}

We are to produce a tree on which the group $\bBS{\m}{\xi}$ acts
transitively. This tree will be constructed from the Bass-Serre
trees of the groups $BS(\m,\xi_n)$. It will be shown that the tree
we construct does not depend on the auxiliary sequence $(\xi_n)_n$.

\begin{figure}[h]
 \begin{center}
   \begin{picture}(250,60)
    \put(7,27){$P$}
    \put(40,30){\circle{40}}
    \put(20,30){\circle*{4}}
    \put(60,32){\vector(0,1){0}}
    \put(67,27){$y$}

    \put(100,35){$G_P = \langle b \rangle \cong \Z$; $G_y = \langle c \rangle \cong \Z$}
    \put(100,15){$c^y = b^{\m}$; $c^{\bar y} = b^{\n}$}
   \end{picture}
 \end{center}
\caption{Baumslag-Solitar groups as graphs of groups} \label{figBS}
\end{figure}
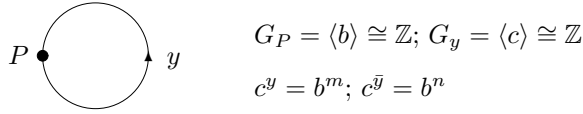

We recall that $BS(\m,\n)$ is the fundamental group of the graph of
groups $(G,Y)$ shown in Figure \ref{figBS} \cite[Section
5.1]{Ser77}.
Notice that $a$ is the element of $\pi_1(G,Y,P)$ associated to the
edge $y$. To be precise, we set the \emph{Bass-Serre tree} of
$BS(\m,\n)$ to be the universal covering associated to $(G,Y)$, the
maximal subtree $P$ and the orientation given by the edge $\bar y$
\cite[Section 5.3]{Ser77}. We choose the edge $\bar y$ instead of
$y$ to minimize the dependence on $\n$ of the set of tree edges.
Denoting by $T$ the Bass-Serre tree of $BS(\m,\n)$, one has
$$
\begin{array}{ccccccc}
                     V(T) & = & BS(\m,\n) / \langle b \rangle & , & \ E_+(T)                  & = & BS(\m,\n) / \langle b^\m \rangle \ , \\
o(w \langle b^\m \rangle) & = & w\langle b \rangle            & , &
t(w \langle b^\m \rangle) & = & wa^{-1}\langle b \rangle \ ,
\end{array}
$$
where $o(e)$ and $t(e)$ denote the origin and terminal vertex of the
edge $e$. The chosen orientation on $T$ is preserved by the
$BS(\m,\n)$-action.

Given $\m \in \Z^*$, $\xi \in \Z_\m$ and $(\xi_n)_n$ a sequence of
rational integers such that $|\xi_n| \to \infty$ and $\xi_n \to \xi$
in $\Z_\m$, we denote by $H_n$ (respectively $H_n^\m$) the subgroup
of $BS(\m,\xi_n)$ generated by $b$ (respectively $b^\m$) and by
$T_n$ the Bass-Serre tree of $BS(\m,\xi_n)$. We set
$$
\begin{array}{ccccc}
Y    & = & \left(\prod\limits_{n\in\N} V(T_n)\right) / \sim   & = & \left(\prod\limits_{n\in\N} BS(\m,\xi_n)/H_n \right) / \sim  \\
Y^\m & = & \left(\prod\limits_{n\in\N} E_+(T_n)\right) / \sim & = &
\left(\prod\limits_{n\in\N} BS(\m,\xi_n)/H_n^\m\right) / \sim
\end{array}
$$
where $\sim$ is defined by $(x_n)_n \sim (y_n)_n \Longleftrightarrow
\exists n_0 \ \forall n \geqslant n_0: \ x_n = y_n$ in both cases.
We now define an oriented graph $X = X_{\m,\xi}$ by
$$
\begin{array}{ccl}
V(X)                     & = & \{ x\in Y : \exists w \in \F \text{ such that } (x_n)_n \sim (wH_n)_n \}  \\
E_+(X)                   & = & \{ y\in Y^\m : \exists w \in \F \text{ such that } (y_n)_n \sim (wH_n^\m)_n \}  \\
o\big( (wH_n^\m)_n \big) & = & (wH_n)_n = \big( o(wH_n^\m) \big)_n  \\
t\big( (wH_n^\m)_n \big) & = & (wa^{-1}H_n)_n =  \big( t(wH_n^\m)
\big)_n
\end{array}
$$
The map $o$ is well defined since $(vH_n^\m)_n \sim (wH_n^\m)_n$
implies $(vH_n)_n \sim (wH_n)_n$. In the other hand, the map $t$ is
well defined since $(vH_n^\m)_n \sim (wH_n^\m)_n$ implies $v^{-1}w
\in H_n^\m$ for $n$ large enough, whence $(va^{-1})^{-1}(wa^{-1}) =
av^{-1}wa^{-1} \in H_n$ for those values of $n$. It follows that
$(va^{-1}H_n)_n \sim (wa^{-1}H_n)_n$. The graph $X$ is thus well
defined and the free group $\F$ acts obviously on it by left
multiplications.
\begin{remark}\label{ultra}
 There is an alternative way to define $X$, which we now describe briefly.
  Consider the $\Z$-trees $V(T_n)$ and define the ultraproduct $V = (\prod V(T_n))/\omega$, which
  is a $^{\omega}\Z$-tree (where $\omega$ is some non-principal ultrafilter over $\N$).
  The ultraproduct group $(\prod BS(m,\xi_n))/\omega$ contains $\bBS{m}{\xi}$ and acts on $V$. Now, the set of vertices of $X$ is equal to subtree of $V$ spanned by the orbit $\bBS{m}{\xi}\cdot v_0$, with $v_0 = [H_n]_\omega$.
   
However, even if this approach would immediately tell us that $X$ is a tree endowed with a $\bBS{m}{\xi}$-action 
(one has only to check that the subtree is a $\Z$-tree), we prefer the more down-to-earth one described above 
for the sake of a self-contained and explicit construction. See \rm \cite{ChGu05, Chi01, Pau04} \it for more information on this second approach.
\end{remark}
It is almost obvious that $X$ is simple (i.e. it has no loop and no
bigon). Indeed a loop (or bigon) in $X$ would immediately provide a
loop (or bigon) in some tree $T_n$. The statement we want to prove
is the following one.

\begin{theorem}\label{TreeAction}
Let $\m \in \Z^*$, $\xi \in \Z_\m$ and $(\xi_n)_n$ a sequence of
rational integers such that $|\xi_n| \to \infty$ and $\xi_n \to \xi$
in $\Z_\m$. The graph $X = X_{\m,\xi}$ (seen here as unoriented)
satisfies the following properties:
\begin{enumerate}
\item[(a)] It is a tree.
\item[(b)] It does not depend (up to equivariant isomorphism) on the choice of the auxiliary sequence
$(\xi_n)_n$.
\item[(c)] The obvious action of the free group $\F$ on $X_{\m,\xi}$ factors through the canonical projection $\F \to \bBS{\m}{\xi}$.
\end{enumerate}
\end{theorem}

Before the proof, we give a simple consequence of \cite[Lemma
6]{Sta06a}, or of Lemma \ref{LemSimCanc} (with $t=h/2$).

\begin{lemma}\label{IdVertices}
Let $(vH_n)_n$ and $(wH_n)_n$ be two vertices of the graph $X$. If
$vH_n = wH_n$ for infinitely many values of $n$, then $(vH_n)_n =
(wH_n)_n$ in $X$.
\end{lemma}

\begin{sproof}{Proof of Theorem \ref{TreeAction}} \textbf{(a)} Let us show first that the graph $X$ is connected. We show by
induction on $|w|$ that any vertex $(wH_n)_n$, with $w \in \F$, is
connected to $(H_n)_n$. The case $|w|=0$ is trivial. If $|w|=\ell >
0$, there exists $x \in \{a^{\pm 1}, b^{\pm 1}\}$ such that $wx$ has
length $\ell - 1$. By induction hypothesis, it is sufficient to show
that $(wH_n)_n$ is connected to $(wxH_n)_n$. If $x=a$, then the edge
$(wxH_n^\m)_n$ connects $(wxH_n)_n$ to $(wH_n)_n$, if $x=a^{-1}$,
then the edge $(wH_n^\m)_n$ connects $(wH_n)_n$ to $(wxH_n)_n$ and
if $x=b^{\pm 1}$, then one has even $(wH_n)_n = (wxH_n)_n$.

\smallskip

Second, we show that $X$ has no closed path without backtracking. We
take a closed path in $X$, with vertices
$$
(v_ 0 H_n)_n, (v_1 H_n)_n, \ldots, (v_\ell H_n)_n=(v_0 H_n)_n \ ,
$$
and we want to prove the existence of some backtracking.

For any $n$, the sequence $v_ 0H_n, v_{1} H_n, \ldots, v_\ell H_n$
defines a path in $T_n$. For $n$ large enough, one has $v_\ell H_n =
v_0 H_n$, so that the path is closed. As the $T_n$'s are trees,
these closed paths all have at least one backtracking. Thus, there
exists $i\in \{0, \ldots, \ell-1\}$ such that $v_{i-1} H_n = v_{i+1}
H_n$ for infinitely many values of $n$ ($v$-indexes are taken modulo
$\ell$). Now, by Lemma \ref{IdVertices}, we obtain $(v_{i-1} H_n)_n
= (v_{i+1} H_n)_n$ in $X$. As $X$ is a simple graph, this means that
the original closed path in $X$ has some backtracking, as desired.

\smallskip

\textbf{(b)} We show now that $X$ does not depend (up to equivariant
isomorphism) on the choice of the sequence $(\xi_n)_n$. Take another
sequence $(\xi'_n)_n$ satisfying both $|\xi'_n| \to \infty$ and
$\xi'_n \to \xi$ in $\Z_\m$ and consider the associated tree $X'$.
We construct the sequence $(\xi''_n)_n$ given by
$$
\xi''_n = \left\{
\begin{array}{cc}
\xi_{\frac{n}{2}} & \text{ if } n \text{ is even } \\
\xi'_{\frac{n-1}{2}} & \text{ if } n \text{ is odd }
\end{array}
\right. \ ,
$$
which satisfies again both $|\xi''_n| \to \infty$ and $\xi''_n \to
\xi$ in $\Z_\m$, and the associated tree $X''$. There are obvious
equivariant surjective graph morphisms $X'' \to X$ and $X'' \to X'$.
We have to show the injectivity of these morphisms, which we can
check on vertices only, for we are dealing with trees. But Lemma
\ref{IdVertices} precisely implies the injectivity on vertices.

\smallskip

\textbf{(c)} Take $w \in \F$ such that $w=1$ in $\bBS{\m}{\xi}$. We
have to prove that $w$ acts trivially on $X$. As $X$ is a simple
graph, we only have to prove that $w$ acts trivially on $V(X)$. Let
$(vH_n)_n$ be a vertex of $X$. For $n$ large enough, we have $w=1$
in $BS(\m,\xi_n)$, so that $wvH_n = vH_n$. Hence we have $w\cdot
(vH_n)_n \sim (vH_n)_n$, as desired.
\end{sproof}

\begin{remark}\label{RemStab}
The action of $\bBS{\m}{\xi}$ on $X$ is transitive and the
stabilizer of the vertex $v_0 = (H_n)_n$ is the subgroup of elements
which are powers of $b$ in all but finitely many $BS(\m,\xi_n)$. It
does not coincide with the subgroup of $\bBS{\m}{\xi}$ generated by
$b$, since the element $ab^\m a^{-1}$ is not in the latter subgroup,
but stabilizes the vertex.
\end{remark}

\begin{remark}\label{RemHNN}
The $\bBS{m}{\xi}$-action on $X$ being transitive on vertices and on
edges, $\bBS{m}{\xi}$ is an HNN-extension with basis
$\operatorname{Stab}(v_0)$. We shall not focus on this structure in
the present article, but leave it for a next paper.
\end{remark}

We end this section by statements about the structure of the tree
$X_{\m,\xi}$. The first one is the analogue of Lemma
\ref{IdVertices} for edges.

\begin{lemma}\label{IdEdges}
Let $(vH_n^\m)_n$ and $(wH_n^\m)_n$ be two edges of the graph $X$.
If one has $vH_n^\m = wH_n^\m$ for infinitely many values of $n$,
then $(vH_n^\m)_n = (wH_n^\m)_n$ in $X$.
\end{lemma}

\begin{proof} By assumption, one has $vH_n = wH_n$ and $va^{-1}H_n = wa^{-1}H_n$ for infinitely many values
of $n$. By Lemma \ref{IdVertices}, we get $(vH_n)_n = (wH_n)_n$ and
$(va^{-1}H_n)_n = (wa^{-1}H_n)_n$. The edges $(vH_n^\m)_n$ and
$(wH_n^\m)_n$ having the same origin and terminal vertex, they are
equal, since $X$ is a simple graph.
\end{proof}

\begin{proposition}\label{TreeStructure}
Let $\m \in \Z^*$ and $\xi \in \Z_\m$. Each vertex of the tree
$X_{\m,\xi}$ has exactly $|\m|$ outgoing edges. More precisely,
(given a sequence $(\xi_n)_n$ of non-zero rational integers such
that $|\xi_n| \to \infty$ and $\xi_n \to \xi$) the edges outgoing
from the vertex $(wH_n)_n$ are exactly $(wH_n^\m)_n, (wbH_n^\m)_n,
\ldots, (wb^{|m|-1}H_n^\m)_n$.
\end{proposition}

\begin{proof} It suffices to treat the case $w=1$. The edges $(H_n^\m)_n$, $(bH_n^\m)_n$,\ldots, $(b^{|\m|-1}H_n^\m)_n$
are clearly outgoing from $(H_n)_n$ and distinct. Let now
$(vH_n^\m)_n$ be an edge outgoing from $(H_n)_n$. In particular, we
have $(vH_n)_n = (H_n)_n$, so that $v=b^{\lambda_n}$ in
$BS(\m,\xi_n)$ for $n$ large enough. There exists necessarily
$\lambda \in \{ 0,\ldots, |\m|-1 \}$ such that we have $\lambda_n
\equiv \lambda$ (mod $\m$) for infinitely many values of $n$, so
that $vH_n^\m = b^{\lambda}H_n^\m$ for infinitely many values of
$n$. By Lemma \ref{IdEdges}, we get $(vH_n^\m)_n = (b^\lambda
H_n^\m)_n$ and we are done.
\end{proof}

\subsection{A structure theorem} \label{SctHaagerup}
\label{Sctstructure}

We denote by $N$ the kernel of the map $q$ appearing in Proposition
\ref{propquot}. We now are able to state
the main results of this section, which are the following:

\begin{theorem}\label{Extension}
Consider the exact sequence (where $N_{\m,\xi}$ is the image of $N$
in $\BS$)
\[
1 \longrightarrow N_{\m,\xi} \longrightarrow \bBS{\m}{\xi}
\overset{q_{\m,\xi}}{\longrightarrow} \Z \wr \Z \longrightarrow 1 \
.
\]
For any $\m\in\Z^*$ and $\xi \in \Z_\m$, the group $N_{\m,\xi} =
\ker q_{\m,\xi}$ is free.
\end{theorem}

\begin{remark}
Since $\Z\wr\Z$ is metabelian, the second derived subgroup of
$\overline{BS}(\m,\xi)$ is a free group. Thus
$\overline{BS}(\m,\xi)$ enjoys the same property as the generalized
Baumslag-Solitar groups \rm\cite[Corollary 2]{Kro90}.
\end{remark}

\begin{sproof}{Proof of Theorem \ref{Extension}} Take $\m\in\Z^*$, $\xi \in \Z_\m$ and $(\xi_n)_n$ a sequence of rational integers such
that $|\xi_n| \to \infty$ and $\xi_n \to \xi$ in $\Z_\m$. Set $X$ to
be the tree constructed in Subsection \ref{Tree}. By \cite[Section
3.3, Theorem 4]{Ser77}, it is sufficient to prove that $N_{\m,\xi}$
acts freely on $X$, i.e. that any $w' \in N$ which stabilizes a
vertex satisfies $w'=1$ in $\BS$.

Let us take $w' \in N$ and $(vH_n)_n$ a vertex of $X$ which is
stabilized by $w'$. Thus $w=v^{-1}w'v$ stabilizes the vertex
$(H_n)_n$, i.e. $w$ is a power of $b$, say $w=b^{\lambda_n}$, in all
but finitely many $BS(\m,\xi_n)$'s. Then the image of $w$ in
$\Gamma(m,\xi_n) = \Z \ltimes_\frac{\n}{\m}
\Z[\frac{\text{gcd}(\m,\xi_n)} {\text{lcm}(\m,\xi_n)}]$ is equal to
$(0,\lambda_n)$ for all but finitely many values of $n$. But, on the
other hand, one has $w\in N$, that is $w=1$ in $\Z\wr\Z$, which
implies $\lambda_n = 0$ for those values of $n$ (since
$\Gamma(m,\xi_n)$'s are marked quotients of $\Z\wr\Z$). Thus, one
has $w=b^0=1$ in all but finitely many $BS(\m,\xi_n)$'s, which gives
$w=1=w'$ in $\BS$.
\end{sproof}

\begin{corollary}\label{Haagerup}
For any $\m\in\Z^*$ and $\xi \in \Z_\m$, the group $\bBS{\m}{\xi}$:
\begin{enumerate}
 \item has the Haagerup property;
 \item is residually solvable.
\end{enumerate}
\end{corollary}

\begin{remark}
The Haagerup property for Baumslag-Solitar groups was already known
\rm\cite{GJ03}\it. It may also be deduced from the fact that the
second derived subgroup is free \rm\cite[Corollary 2]{Kro90}.
\end{remark}

\begin{sproof}{Proof of Corollary \ref{Haagerup}} Looking at the exact sequence
$$
1 \longrightarrow N_{\m,\xi} \longrightarrow \bBS{\m}{\xi}
\overset{q}{\longrightarrow} \Z \wr \Z \longrightarrow 1 \ ,
$$
we see that the quotient group is amenable (it is even metabelian)
and the kernel group has the Haagerup property by Theorem
\ref{Extension}. By \cite[Example 6.1.6]{CCJ+01}, we conclude that
$\BS$ has the Haagerup property. As a free group is residually
solvable, $\bBS{\m}{\xi}$ is then the extension of a residually
solvable group by a solvable one and hence is residually solvable by
\cite[Lemma 1.5]{Gru57}.
\end{sproof}

\section{Presentations of the limits}\label{SctPres1}

The purpose of this Section is, first, to prove that the $\BS$'s are
almost never finitely presented (Theorem \ref{infPres}), and second,
to give presentations, whose construction is based on actions on
trees described in Section \ref{Tree} (Theorem \ref{ThmPres1}).

\subsection{Most of the limits are not finitely presented}

Our first goal in this Section is to prove:

\begin{theorem}\label{infPres}
For any $\m \in \Z^*$ and $\xi \in \Z_\m \setminus \m\Z_\m$, the
group $\ov{BS}(\m,\xi)$ is not finitely presented.
\end{theorem}

Notice that the Theorem excludes also the existence of a finite
presentation of such a group with another generating set. See for
instance \cite[Proposition V.2]{Har00}. By Theorem \ref{ThmBSm},
there is only one remaining case, the case $\xi=0$, where it is
still unknown whether $\overline{BS}(\m,0)$ is finitely presented or
not.
 We nevertheless make the following remark.
\begin{remark}\label{RemUniqZ1}
For $|\m|=1$, the limits $\bBS{\pm 1}{\xi}$ are not finitely
presented.
\end{remark}

\begin{proof}
The result \cite[Theorem 2]{Sta06a} implies $\bBS{\pm 1}{\xi} =
\Z\wr\Z$ for the unique element $\xi\in\Z_{\pm 1}$ and the Baumslag's
result \cite{Bau61} on the presentations of wreath products
ensures that $\Z\wr\Z$ is not finitely presented.
\end{proof}

\begin{lemma}\label{lemInfPres}
Let $\m \in \Z^*$ and $\xi \in \Z_\m \setminus \m\Z_\m$. Let $\ell$
be the maximal exponent in the decomposition of $\m$ in prime
factors and set $d= \text{gcd}(\m,\xi)$, $\m_1 = \m/d$.
\begin{enumerate}
  \item[(a)] There exists a sequence $(\xi_n)_n$ in $\Z^*$ such that for all $n \geqslant 1$ one has
  $|\xi_n|>|\xi_{n-1}|$ and
    $$\begin{array}{cl}
        \xi_n \equiv \xi & (\text{mod } \m^n\Z_\m) \ ; \\
        \xi_n \not\equiv \xi & (\text{mod } \m_1^{\ell n + 1}d\Z_\m) \ .
    \end{array}$$

  \item[(b)] This sequence satisfies $|\xi_n|\to\infty$, $\xi_n \to \xi$ and
    $$\begin{array}{cll}
        \xi_n \equiv \xi_r & (\text{mod } \m^n) & \forall r\geqslant n \ ; \\
        \xi_n \not\equiv \xi_{\ell n+1} & (\text{mod } \m_1^{\ell n + 1}d) & \forall n \ .
    \end{array}$$
\end{enumerate}

\end{lemma}

\begin{proof}
\textbf{(a)} Let $p$ be a prime factor of $\m_1$ (there exists one,
for $\xi \not \in \m\Z_\m$). The sequence $(\xi_n)_n$ is constructed
inductively. We choose for $\xi_0$ any non-zero rational integer
such that $\xi_0 - \xi \not \in \m\Z_\m$. At the $n$-th step, we
begin by noticing that the exponent of $p$ in the decomposition of
$\m^n$ (respectively $m_1^{\ell n +1}d$) is at most $\ell n$
(respectively at least $\ell n + 1$). Hence, $m_1^{\ell n +1}d$ is
not a divisor of $\m^n$, so that there exists $\alpha\in\Z$ with
$\xi \equiv \alpha \ (\text{mod } \m^n \Z_m)$ but $\xi \not\equiv
\alpha \ (\text{mod } \m_1^{\ell n + 1}d \Z_m)$.

Notice now that we may replace $\alpha$ by any element of the class
$\alpha + \m^n \m_1^{\ell n + 1}d\Z$, so that it suffices to choose
$\xi_n$ among the elements $\beta$ in the latter class which satisfy
$|\beta|
> |\xi_{n-1}|.$

\textbf{(b)} The properties $\xi_n \equiv \xi \ (\text{mod }
\m^n\Z_\m)$ and $|\xi_n|>|\xi_{n-1}|$ imply clearly $\xi_n \to \xi$,
$|\xi_n|\to\infty$ and $\xi_n \equiv \xi_r \ (\text{mod } \m^n)$ for
$r\geqslant n$ (for the latter one, Proposition \ref{mAdic} (d) is
used).

Finally, combining the properties $\xi_{\ell n +1} \equiv \xi \
(\text{mod } \m^{\ell n+1}\Z_\m)$ and $\xi_n \not\equiv \xi \
(\text{mod } \m_1^{\ell n + 1}d\Z_\m)$ gives $\xi_n \not\equiv
\xi_{\ell n+1} \ (\text{mod } \m_1^{\ell n +1}d)$.
\end{proof}

\medskip

\begin{sproof}{Proof of Theorem \ref{infPres}}
The hypothesis $\xi \in \Z_\m \setminus \m\Z_\m$ implies $|\m|
\geqslant 2$. Take $\ell, d, \m_1$ and a sequence $(\xi_n)_n$ as in
Lemma \ref{lemInfPres}. One has then $BS(\m,\xi_n) \to \BS$. It is
thus sufficient by Lemma \ref{presFinQuot} to prove that the
$BS(\m,\xi_n)$'s are not marked quotients of $\BS$ (for $n$ large
enough).

Notice now that (for $n$ large enough) one has $\text{gcd}(\m,\xi_n)
= d$ since $\xi_n \equiv \xi \ (\text{mod } \m\Z_\m)$ holds. Then,
for $n\geqslant 1$, we define the words
$$
w_n = a^{n+1}b^\m a^{-1}b^{-\xi_n}a^{-n}b
a^{n+1}b^{-\m}a^{-1}b^{\xi_n}a^{-n}b^{-1}.
$$
By Lemma 3 of \cite{Sta06a}, we have $w_n=1$ in $BS(m,k)$ if and
only if $k \equiv \xi_n \,(\text{mod }m_1^nd)$. By Part (b) of Lemma
\ref{lemInfPres}, we have then $w_n = 1$ in $BS(\m,\xi_r)$ for all
$r\geqslant n$, hence $w_n = 1$ in $\BS$ (for $n$ large enough). On
the other hand we get $w_{\ell n+1} \neq 1$ in $BS(\m,\xi_n)$ the
same way, so that $BS(\m,\xi_n)$ is not a marked quotient of $\BS$
(for $n$ large enough).
\end{sproof}

\subsection{Presentations of the limits}
As we noticed before (see Remark \ref{RemHNN}), the limits
$\bBS{\m}{\xi}$ are HNN-extensions. This would give us presentations
of the limits, provided that we would have presentations for the
vertex and edges stabilizers. Nevertheless, rather than to follow
this approach and seek presentations for the stabilizers, we
establish directly a presentation for the whole group.

For $\m\in\Z^*$ and $\xi \in \Z_\m$, we define the set ${\cal
R}={\cal R}_{\m,\xi}$ by
\begin{eqnarray*}
{\cal R}_{\m,\xi} & = & \big\{ w\bar w : w = ab^{\alpha_1} \cdots
ab^{\alpha_k} a^{-1}b^{\alpha_{k+1}} \cdots
a^{-1}b^{\alpha_{2k}} \\
 & & \ \text{ with } k\in \N^*, \alpha_i \in \Z \ (i=1, \ldots, 2k) \text{ and } w\cdot v_0 = v_0 \big\}
\end{eqnarray*}
where $v_0$ is the favoured vertex $(H_n)_n$ of the tree
$X_{\m,\xi}$ defined in Section \ref{Tree}.

Recall that the stabilizer of the vertex $v_0$ consists of elements
which are powers of $b$ in all but finitely many $BS(\m,\xi_n)$,
where $(\xi_n)_n$ is any sequence of rational integers such that
$|\xi_n| \to \infty$ and $\xi_n \to \xi$ in $\Z_\m$ (see Remark
\ref{RemStab}). It follows that we have $w\cdot v_0 = v_0
\Leftrightarrow \bar w\cdot v_0 = v_0$. The aim is now to prove the
following result.
\begin{theorem}\label{ThmPres1}
For all $\m \in \Z^*$ and $\xi \in \Z_\m$, the marked group
$\bBS{\m}{\xi}$ admits the presentation $\Pres{a,b}{{\cal
R}_{\m,\xi}}$.
\end{theorem}
Set $\Gamma = \Pres{a,b}{{\cal R}_{\m,\xi}}$ until the end of this
Section. The elements of ${\cal R}_{\m,\xi}$ are
trivial in $\bBS{\m}{\xi}$: indeed, if $w$ fixes the favoured vertex $H_n$ in
the Bass-Serre tree of $BS(m,\xi_n)$, then $w\bar w$ is trivial in
$BS(m,\xi_n)$. This gives a marked (hence surjective) homomorphism $\Gamma \to
\bBS{\m}{\xi}$. Theorem \ref{ThmPres1} is thus reduced to the following
proposition, which gives the injectivity.
\begin{proposition}\label{PropPres1}
Let $w$ be a word on the alphabet $\{ a,a^{-1},b,b^{-1} \}$. If one
has $w = 1$ in $\bBS{\m}{\xi}$, then the equality $w = 1$ also holds
in $\Gamma$.
\end{proposition}

Before proving Proposition \ref{PropPres1}, we need to introduce
some notions which admit geometric interpretations.

\medskip

\textbf{From words to paths.} Let us call \emph{path} (in a graph)
any finite sequence of vertices such that each of them is adjacent
to the preceding one. Let $w$ be any word on the alphabet $\{
a,a^{-1},b,b^{-1} \}$. It defines canonically a path in the Cayley
graph of $\Gamma$ (or $\bBS{\m}{\xi}$) which starts at the trivial
vertex. Let us denote those paths by $p_\Gamma(w)$ and
$p_{\BBS}(w)$. The map $\Gamma \to \bBS{\m}{\xi}$ defines a graph
morphism which sends the path $p_\Gamma(w)$ onto the path
$p_{\BBS}(w)$.

The word $w$ defines the same way a finite sequence of vertices in
$X_{\m,\xi}$ starting at $v_0$ and such that each of them is equal
or adjacent to the preceding one. Indeed, let $f$ be the map $V
\big( Cay(\bBS{\m}{\xi},(a,b)) \big) \to V(X_{\m,\xi})$ defined by
$f(g) = g\cdot v_0$ for any $g \in \bBS{\m}{\xi}$. If $g,g'$ are
adjacent vertices in $Cay(\bBS{\m}{\xi},(a,b))$, then $f(g)$ and
$f(g')$ are either adjacent (case $g'=ga^{\pm 1}$), or equal (case
$g'=gb^{\pm 1}$). The sequence associated to $w$ is the image by $f$
of the path $p_{\BBS}(w)$. Now, deleting consecutive repetitions in
this sequence, we obtain a path that we denote by $p_X(w)$.

It follows that if the word $w$ satisfies $w = 1$ in $\bBS{\m}{\xi}$
(or, stronger, $w=1$ in $\Gamma$), then the path $p_X(w)$ is closed
(i.e. its last vertex is $v_0$).

\medskip

\textbf{Height and Valleys.} Recall that one has a homomorphism
$\sigma_a$ from $\bBS{\m}{\xi}$ onto $\Z$ given by $\sigma_a(a)=1$
and $\sigma_a(b)=0$. Given a vertex $v$ in $X_{\m,\xi}$, we call
\emph{height} of $v$ the number $h(v) = \sigma_a(g)$ where $g$ is
any element of $\bBS{\m}{\xi}$ such that $g\cdot v_0 = v$. It is
easy to check that any element $g'$ of $\bBS{\m}{\xi}$ which defines
an elliptic automorphism of $X_{\m,\xi}$ satisfies $\sigma_a(g') =
0$, so that the height function is well-defined. It is clear from
construction that the height difference between two adjacent
vertices is $1$.

Given $L \geqslant 1$ and $k\geqslant 1$, we call
\emph{$(L,k)$-valley} any path $p$ in $X_{\m,\xi}$ such that one
has:
\begin{itemize}
    \item $p = (v_0,v_1, \ldots, v_L = \nu_0, \nu_1, \ldots, \nu_{2k})$, where $v_0$ is the favoured
    vertex;
    \item $h(v_0) = 0 = h(\nu_k)$ and $h(\nu_0) = -k = h(\nu_{2k})$;
    \item $h(v)<0$ for any other vertex $v$ of $p$;
    \item $\nu_0 = \nu_{2k}$.
\end{itemize}

Given a $(L,k)$-valley $p = (v_0,v_1, \ldots, v_L = \nu_0, \nu_1,
\ldots, \nu_{2k})$, the subpaths $(\nu_0, \ldots, \nu_k)$ and
$(\nu_k, \ldots, \nu_{2k} = \nu_0)$ have to be geodesic, for the
height difference between $\nu_0 = \nu_{2k}$ and $\nu_k$ is $k$.
Thus, one has $\nu_1 = \nu_{2k-1}, \ldots, \nu_{k-1} = \nu_{k+1}$.

\begin{lemma}\label{LemPres1}
Let $w$ be a word on the alphabet $\{ a,a^{-1},b,b^{-1} \}$ such
that the path $p_X(w)$ is a $(L,k)$-valley, say $p_X(w) = (v_0,v_1,
\ldots, v_L = \nu_0, \nu_1, \ldots, \nu_{2k})$. There exists a word
$w'$ such that the equality $w'=w$ holds in $\Gamma$ and the path
$p_X(w')$ is $(v_0,v_1, \ldots, v_L)$.
\end{lemma}

\begin{proof} We argue by induction on $L$.

\textbf{Case $L=1$:} In that case, one has $k=1$, $p_X(w) = (v_0,
v_1=\nu_0, \nu_1, \nu_2)$. Up to replacing $w$ by a word which
defines the same element in $\F$ (hence in $\Gamma$) and the same
path in $X$, we may assume to have $w = b^{\alpha_0} a^{-1}
b^{\alpha_1} a b^{\beta_1} a^{-1} b^{\beta_2}$. Set $r= a
b^{-\beta_1} a^{-1} b^{-\alpha_1} a b^{\beta_1} a^{-1}
b^{\alpha_1}$. Since $\nu_0$ and $\nu_2$ are equal, the subword $a
b^{\beta_1} a^{-1}$ (of $w$) defines a closed subpath in $X$, so
that we obtain $a b^{-\beta_1} a^{-1} b^{-\alpha_1} \cdot v_0 =
v_0$. Consequently, we get $r \in {\cal R}$, whence $r=1$ in
$\Gamma$. Inserting $r$ in next to last position, we obtain
$$
w \underset{\Gamma}{=} b^{\alpha_0 + \beta_1} a^{-1} b^{\alpha_1 +
\beta_2} = : w' \ .
$$
This equality also implies that the paths $p_X(w)$ and $p_X(w')$
have the same endpoint. Hence one has $p_X(w') = (v_0,\nu_2) =
(v_0,v_1)$ and we are done.

\textbf{Induction step:} We assume $L>1$ to hold. Up to replacing
$w$ by a word which defines the same element in $\F$ (hence in
$\Gamma$) and the same path in $X$, we may write
$$
w = b^{\alpha_0} a^{\varepsilon_1} b^{\alpha_1} \cdots
a^{\varepsilon_L} b^{\alpha_L}\cdot a b^{\beta_1} \cdots a
b^{\beta_k} a^{-1} b^{\beta_{k+1}} \cdots a^{-1} b^{\beta_{2k}}
$$
with $\varepsilon_i = \pm 1$ and $\alpha_i\in \Z$ for all $i$. We
distinguish two cases:
\begin{enumerate}
    \item[(1)] the vertex $v_{L-1}$ is higher than $v_L$ (i.e. $\varepsilon_L = -1$);
    \item[(2)] the vertex $v_{L-1}$ is lower than $v_L$ (i.e. $\varepsilon_L = 1$).
\end{enumerate}

\textit{Case (1):} Set
\begin{eqnarray*}
z = a b^{-\beta_{2k-1}} \cdots a b^{-\beta_k} a^{-1}
b^{-\beta_{k-1}} \cdots a^{-1} b^{-\beta_1}
a^{-1} b^{-\alpha_L} & \text{ and } & r = z \bar z \ .
\end{eqnarray*}
Since $\nu_0$ and $\nu_{2k}$ are equal, the subword $a b^{\beta_1}
\cdots a b^{\beta_k} a^{-1} b^{\beta_{k+1}} \cdots a^{-1}
b^{\beta_{2k-1}} a^{-1}$ (of $w$) defines a closed subpath in $X$,
so that (when considered as a word on its own right) it stabilizes
the vertex $v_0$. Inverting it, we get
$$
a b^{-\beta_{2k-1}} \cdots a b^{-\beta_k} a^{-1} b^{-\beta_{k-1}}
\cdots a^{-1} b^{-\beta_1} a^{-1} \cdot v_0 = v_0 \ .
$$
It implies $r \in {\cal R}$, whence $r=1$ in $\Gamma$. Inserting $r$
in next to last position, we obtain
\begin{eqnarray*}
w \underset{\Gamma}{=} w^* & := & b^{\alpha_0} a^{\varepsilon_1}
b^{\alpha_1} \cdots a^{\varepsilon_{L-1}}
b^{\alpha_{L-1}+\beta_{2k-1}} \cdot \\
 & &  a b^{\beta_{2k-2}} \cdots a b^{\beta_k} a^{-1} b^{\beta_{k-1}} \cdots a^{-1} b^{\beta_1} \cdot \\
 & &  a^{-1} b^{\alpha_L+\beta_{2k}} \ .
\end{eqnarray*}
We write $w^*=w''a^{-1} b^{\alpha_L+\beta_{2k}}$. Since the words
$w$ and $w^*$ begin the same way and since one has $w = w^*$ in
$\Gamma$, the path $p_X(w^*)$ has the form
$$
(v_0,v_1, \ldots, v_{L-1} = \omega_0, \omega_1, \ldots,
\omega_{2k-2}, v_L) \ .
$$
Since $v_{L-1}$ is higher than $v_L$, we have $h(v_{L-1}) = -(k-1)$.
Contemplating $w^*$, one sees that we have $h(\omega_{k-1}) = 0$,
$h(\omega_{2k-2}) = -(k-1)$, so that the subpaths $(\omega_0,
\ldots, \omega_{k-1})$ and $(\omega_{k-1}, \ldots, \omega_{2k-2},
v_L)$ are geodesic. On the other hand, the geodesic between
$\omega_{k-1}$ and $v_L$ passes through $v_{L-1}=\omega_0$, so that
we have $\omega_{2k-2}=v_{L-1}$. It follows that $p_X(w'')$ has the
form $(v_0,v_1, \ldots, v_{L-1} = \omega_0, \omega_1, \ldots,
\omega_{2k-2} = v_{L-1})$, so that it is a $(L-1,k-1)$-valley. We
apply the induction hypothesis to $w''$ and get a word $w'''$ such
that $w''=w'''$ holds in $\Gamma$ and $p_X(w''') = (v_0, \ldots,
v_{L-1})$. It suffices to set $w'=w'''a^{-1}
b^{\alpha_L+\beta_{2k}}$ to conclude.

\textit{Case (2):} In this case, we have $h(v_L) = -k$ and
$h(v_{L-1}) = -(k+1)$, so that the edge linking these vertices goes
from $v_L$ to $v_{L-1}$. By Proposition \ref{TreeStructure}, there
exists $\lambda \in \{ 0, \ldots, |\m|-1 \}$ such that $wb^\lambda
a^{-1}\cdot v_0 = v_{L-1}$. Set $w^* = wb^\lambda a^{-1}$ and $\tilde w=
w^*ab^{-\lambda}$, so that $w=\tilde w$ holds in $\Gamma$. The path
$p_X(w^*)$ is an $(L-1,k+1)$-valley, so that we may apply the
induction hypothesis to $w^*$. This gives a word $w''$ such that
$w''=w^*$ in $\Gamma$ and $p_X(w'') = (v_0, \ldots, v_{L-1})$. We
conclude by setting $w'=w''ab^{-\lambda}$.
\end{proof}

\medskip

\begin{sproof}{Proof of Proposition \ref{PropPres1}} Let $w$ be a word which
defines the trivial element in $\bBS{\m}{\xi}$. The path $p_X(w)$ is closed;
we denote its length by $\ell$, which is equal to the number of letters $a^{\pm 1}$
occuring in $w$. We argue by induction on $\ell$.

\textbf{Case $\ell = 0$:} In this case, $w$ has no $a^{\pm 1}$ letter; thus,
one has $w=b^{\alpha_0}$ in $\F$, so that $b^{\alpha_0} = 1$ in
$\bBS{\m}{\xi}$. It follows $\alpha_0 = 0$, hence $w=1$ in $\Gamma$.

\textbf{Induction step ($\ell>0$):} Let us write $p_X(w) = (v_0, v_1, \ldots,
v_{\ell -1}, v_\ell = v_0)$. If $w$ is not freely cyclically reduced, let $w_r$
denote the word obtained by freely cyclically reducing $w$. Then, either some
letter $a$ or $a^{-1}$ has been deleted, in which case $p_X(w_r)$ is strictly
shorter than $p_X(w)$ and we may apply the induction hypothesis to get
$w=w_r=1$ in $\Gamma$, or $p_X(w_r)$ has length $\ell$ and it is equivalent
to deal with $w_r$ instead of $w$.

Hence, we may assume $w$ to be freely cyclically reduced and write
\[
w = b^{\alpha_0} a^{\varepsilon_1} b^{\alpha_1} \cdots
a^{\varepsilon_\ell} b^{\alpha_\ell}
\]
with $\varepsilon_i = \pm 1$ and $\alpha_i \in \Z$ for all $i$.
We may moreover assume $h(v_i) \leqslant 0$ for all $i$: indeed, if the vertex
$v_k$ has maximal height (among the $v_i$'s), then it suffices to replace $w$
by the free reduction of the conjugate $x^{-1}wx$, where $x$ is the prefix
$b^{\alpha_0} \cdots  b^{\alpha_{k-1}}$ of $w$. Note that this operation
preserves the length of $p_X(w)$ and the fact that $w$ is freely cyclically
reduced.

Denote by $k_0=0 < k_1 < \ldots < k_s=\ell$ the indices $k'$ such that $h(v_{k'}) = 0$.
We now distinguish two possibilities:
\begin{enumerate}
    \item[(1)] The path $p_X(w)$ backtracks at some $v_{k_i}$ (i.e. $\exists i$ with $1\leqslant i \leqslant s-1$
    such that $v_{k_i+1} = v_{k_i-1}$.)
    \item[(2)] The path $p_X(w)$ does not backtrack at any $v_{k_i}$ (i.e. $\forall i$ with $1\leqslant i \leqslant s-1$
    one has $v_{k_i+1} \neq v_{k_i-1}$.)
\end{enumerate}
\textit{Case (1):} Let $i$ be an index (with $1\leqslant i \leqslant
s-1$) such that $v_{k_i+1} = v_{k_i-1}$. The subword $w' =
a^{\varepsilon_{k_{i-1}+1}} b^{\alpha_{k_{i-1}+1}} \cdots
a^{\varepsilon_{k_i}} b^{\alpha_{k_i}}a^{\varepsilon_{k_i+1}}$
defines by construction a $(k_i - k_{i-1} -1,1)$-valley in the tree
$X$. Lemma \ref{LemPres1} furnishes then a word $w''$ such that
$w''=w'$ holds in $\Gamma$ and the path $p_X(w'')$ is strictly
shorter than $p_X(w')$. We construct a word $w^*$ by replacing $w'$
by $w''$ in $w$. The path $p_X(w^*)$ is strictly shorter than
$p_X(w)$ and one has $w^*=w$ in $\Gamma$. Applying the induction
hypothesis to $w^*$, we get $w^* = 1$ in $\Gamma$, hence $w=1$ in
$\Gamma$.

\textit{Case (2):} Let us recall that for all $n\in \Z^*$, the following
diagram of marked morphisms is commutative (see Section \ref{NotConv} and
Proposition \ref{propquot}).
$$
\xymatrix{
\F \ar@{->>}[d] \ar@{->>}[r] & \Gamma \ar@{->>}[r] & \overline{BS}(\m,\xi)
\ar@{->>}[r] & \Z\wr\Z \ar@{->>}[d] \\
BS(\m,\n)  \ar@{->>}[rrr] & & & \Gamma(\m,\n) = \Z \ltimes_{\frac{\n}{\m}}
\Z[\frac{\text{gcd}(\m,\n)}{\text{lcm}(\m,\n)}]}
$$
The image of $w$ in $\Z\wr\Z$ is $(\sum_{i=1}^\ell  \varepsilon_i, \sum_{i=0}^\ell
\alpha_i t^{h(v_i)})$. Since $w=1$ in $\overline{BS}(\m,\xi)$, this implies
$\sum_{i=1}^\ell  \varepsilon_i = 0$ and $\sum_{i=0}^\ell \alpha_i t^{h(v_i)} = 0$.

Suppose first (by contradiction) that $s=1$. All vertices of
$p_X(w)$ but $v_0$ and $v_\ell$ have strictly negative height. The image of $w$
in $\Z\wr\Z$ has then the form $(0, \alpha_0+\alpha_\ell + \sum_{j<0}
\beta_j t^j)$, so that $\alpha_0 = -\alpha_\ell$. Since $h(v_1) = -1 =
h(v_{\ell - 1})$, we also have $\varepsilon_1 = -1$ and
$\varepsilon_\ell = 1$. We thus see that $w$ is not freely cyclically reduced,
a contradiction.

Hence, we have $s\geqslant 2$. There exists some $i$ in $\{0,\ldots,
s-1 \}$ such that one has $v_{k_i} = v_{k_{i+1}}$: either one has $v_{k_i} =
v_0$ for all $i$, or we take $i$ such that $v_{k_i}$ is farthest from $v_0$.
We set
$$
w^* := a^{\varepsilon_{k_i+1}} b^{\alpha_{k_i+1}} \cdots
a^{\varepsilon_{k_{i+1}-1}} b^{\alpha_{k_{i+1}-1}}
a^{\varepsilon_{k_{i+1}}} \ .
$$
It is a subword of $w$, not containing all $a^{\pm 1}$ letters since
$s\geqslant 2$. Moreover, the path $p_X(w^*)$ has the form
$(v_0^*=v_0, v_1^*, \ldots, v_{\ell^*-1}^*, v_{\ell^*}^* = v_0)$.
Let us now fix a sequence $(\xi_n)_n$ of non-zero rational integers
such that $|\xi_n| \to \infty$ and $\xi_n \to \xi$ in $\Z_\m$. By
Theorem \ref{TreeAction}, we may assume the tree $X$ to be
constructed from the sequence $(\xi_n)_n$. The word $w^*$
stabilizing $v_0$, we get $w^* = b^{\lambda_n}$ in $BS(\m,\xi_n)$
for some $\lambda_n \in \Z$, which implies $w^* = b^{\lambda_n}$ in
$\Gamma(\m,\xi_n)$, for $n$ large enough. Now, the above diagram
implies
$$
\lambda_n = \sum\limits_{j=1}^{\ell^*-1} \alpha_{k_i+j} \left(
\frac{\xi_n}{\m} \right)^ {h(v_j^*)}
$$
for those values of $n$. Then, taking absolute values, this gives
$$
|\lambda_n| \leqslant \sum\limits_{j=1}^{\ell^*-1} |\alpha_{k_i+j}
\left( \frac{\xi_n}{\m} \right)^ {h(v_j^*)}| \leqslant
\frac{|\m|}{|\xi_n|}\sum\limits_{j=1}^{\ell^*-1} |\alpha_{k_i+j}| \ .
$$
It follows that $|\lambda_n| < 1$ holds for $n$ large enough, since
$|\xi_n|$ tends to $\infty$. For those values of $n$, we get $w^* =
b^0 = 1$ in $BS(\m,\xi_n)$. Consequently, we get $w^* = 1$ in
$\bBS{\m}{\xi}$.

Since $w^*$ does not contain all $a^{\pm 1}$ letters of $w$, the
path $p_X(w^*)$ is strictly shorter than $p_X(w)$, so that we apply
the induction hypothesis to $w^*$ and get $w^* = 1$ in $\Gamma$.
Erasing $w^*$ in $w$, we get
$$
w \underset{\Gamma}{=} w' = b^{\alpha_0} a^{\varepsilon_1}
b^{\alpha_1} \cdots a^{\varepsilon_{k_i}}
b^{\alpha_{k_i}+\alpha_{k_{i+1}}} a^{\varepsilon_{k_{i+1}+1}}
b^{\alpha_{k_{i+1}+1}} \cdots a^{\varepsilon_\ell} b^{\alpha_\ell} \
.
$$
Applying the induction hypothesis to $w'$, we get $w=1$ in $\Gamma$.
\end{sproof}

\bibliographystyle{alpha}
\newcommand{\etalchar}[1]{$^{#1}$}
\def\cprime{$'$} \def\cprime{$'$} \def\cprime{$'$} \def\cprime{$'$}

\medskip

Authors addresses:

\medskip

L. G. G\"ottingen Universit\"at, Mathematisches Institut, Bunsenstrasse 3-5, 
37073 G\"ottingen, Germany, guyot@uni-math.gwdg.de

\medskip

Y. S. Laboratoire de Math\'ematiques, Universit\'e Blaise Pascal,
Campus Universitaire des C\'ezeaux, F-63177 Aubi\`ere cedex, France,
yves.stalder@math.univ-bpclermont.fr

\end{document}